\documentclass[12pt,reqno]{amsart}
\usepackage{amsthm}
\usepackage{amstext}
\usepackage{amssymb}
\usepackage{mathrsfs}
\usepackage{mathabx}
\usepackage[english]{babel}
\usepackage{enumerate}
\usepackage{xfrac}
\usepackage{mathabx}
\usepackage{mathtools,textcomp}
\usepackage{xcolor}
\usepackage[all]{xy}
\usepackage{todonotes}
\usepackage[T1]{fontenc}

\theoremstyle{plain}
\newtheorem{theorem}{Theorem}[section]

\newtheorem{proposition}[theorem]{Proposition}

\newtheorem{lemma}[theorem]{Lemma}

\theoremstyle{definition}
\newtheorem{definition}[theorem]{Definition}

\newtheorem{example}[theorem]{Example}
\newtheorem{remark}[theorem]{Remark}

\theoremstyle{remark}

\numberwithin{equation}{section}

\newcommand{\supp}{\mathrm{supp}}
\newcommand{\set}[2]{\bigl\{#1\bigm|#2\bigr\}}

\newcommand{\seqx}[2]{\left( #1_#2\right)_{#2\in\IN}}
\newcommand{\En}{\mathscr{E}}
\newcommand{\F}{\mathscr{F}}

\newcommand{\efloc}{\mathscr{F}_{\mathrm{loc}}}

\newcommand{\pt}[1]{{#1}^{\bullet}}

\newcommand{\B}{\mathsf{B}}

\newcommand{\IN}{\mathbb{N}}
\newcommand{\IR}{\mathbb{R}}
\newcommand{\m}{\mathfrak{m}}
\newcommand{\loc}{\mathrm{loc}}
\newcommand{\Id}{\mathrm{d}}
\newcommand{\dm}{\mathsf{d}}
\newcommand{\IP}{\mathbb{P}}
\newcommand{\IWW}{\mathscr{W}}
\newcommand{\dom}{\mathrm{Dom}}

\title[Subharmonicity on locally smoothing spaces]{A new notion of subharmonicity on locally smoothing spaces, and a conjecture by Braverman, Milatovic, Shubin }%

\author{Batu Güneysu}

\author{Stefano Pigola}

\author{Peter Stollmann}

\author{Giona Veronelli}

\begin{document}

\begin{abstract} Given a strongly local Dirichlet space and $\lambda\geq 0$, we introduce a new notion of $\lambda$--subharmonicity for $L^1_\loc$--functions, which we call \emph{local $\lambda$--shift defectivity}, and which turns out to be equivalent to distributional $\lambda$--subharmonicity in the Riemannian case. We study the regularity of these functions on a new class of strongly local Dirichlet, so called locally smoothing spaces, which includes Riemannian manifolds (without any curvature assumptions), finite dimensional RCD spaces, Carnot groups, and Sierpinski gaskets. As a byproduct of this regularity theory, we obtain in this general framework a proof of a conjecture by Braverman, Milatovic, Shubin on the positivity of distributional $L^q$-solutions of $\Delta f\leq f$ for complete Riemannian manifolds.
\end{abstract}

\maketitle

\section{Introduction}

Given $q\in [1,\infty]$ one says that a connected Riemannian manifold $X$ is \emph{$L^q$-positivity preserving}, if for every $f\in L^q(X)$ the following implication holds true:
\begin{align}\label{ppp}
\text{$\Delta f\leq f$ in the sense of distributions $\Rightarrow$  $f$ is nonnegative. }
\end{align}
Here $\Delta=\sum_{ij}g^{ij}\partial_i\partial_j$ is the negative-definite Laplace-Beltrami operator. The importance and subtlety of this property is reflected at least by the following observations:
\begin{itemize}
\item if $X$ is $L^2$-positivity preserving, then $-\Delta$ (defined on $C^\infty_c(X)$) is essentially self-adjoint (in $L^2(X))$ \cite{bms},
\item $X$ is stochastically complete, if and only if $X$ is $L^\infty$-positivity preserving \cite{bismar},
\item there exist complete $X$'s which are not $L^1$-positivity preserving \cite{bismar}, 
\item a conjecture by Braverman-Milatovic-Shubin (BMS) from 2002 \cite{bms} states that if $X$ is complete, then $X$ is $L^2$-positivity preserving.
\end{itemize}
The BMS-conjecture was formulated in an $L^2$-setting as the authors were interested in essential self-adjointness problems and was based on the classical fact that $-\Delta$ is essentially self-adjoint, if $X$ is complete. More generally, one can ask, whether completeness implies $L^q$-positivity preservation for all $q\in (1,\infty)$.\vspace{1mm}

Let us point out that the main technical problem in proving such a conjecture is the possible lack of regularity of $f$ above. For if $f$ was in $W^{1,2}_\loc(X)$, then one can integrate by parts once to see that the inequality in (\ref{ppp}) is equivalent to 
\begin{align}\label{ppp2}
\text{$\int_X (\nabla f ,\nabla  \phi)+\int_X f\phi \geq 0 \;$ for all $0\leq \phi\in  W^{1,2}_c(X)$, }
\end{align}
in other words, $f$ is weakly $1$-superharmonic. If $f$ was even in $W^{1,2}_0(X)$, then it would be easy to show that (\ref{ppp2}) implies $f\ge 0$ using first order cut-off functions in the sense of \cite{batu2}, the existence of which is equivalent to the completeness of $X$. \vspace{1mm}

On another note, it is straightforward to see that, if $X$ admits a sequence of Laplacian cut-off functions in the sense of \cite{batu2}, then $X$ is $L^q$-positivity preserving for all $q\in (1,\infty)$. However, in general, this procedure (which relies on two integrations by parts \cite{batu1,batu2}) requires \cite{bianchi} a lower control of the Ricci curvature in addition to the completeness of $X$. \vspace{1mm}

It is well-known that a convenient abstract setting to formulate an inequality of the form (\ref{ppp2}) is provided by a (regular and symmetric) Dirichlet space $(X,\m,\mathscr{E}, \F(X))$ \cite{fukushima_classic,maroeckner,fukushima,chenfukushima}. Indeed, given a nonnegative real number $\lambda$, one has a natural notion of weakly $\lambda$-subharmonic and weakly $\lambda$-superharmonic functions at hand \cite{sturm1}, where the role of $W^{1,2}_\loc(X)$ and $W^{1,2}_c(X)$ is, respectively, played by $\F_\loc(X)$ and $\F_c(X)$, with $\F(X)$ the domain of definition of the underlying quadratic form $\mathscr{E}$. Note that in the Riemannian case, one has
$$
\F(X)=W^{1,2}_0(X),\quad \En(f,g)=\int_X (\nabla f,\nabla g),
$$
the usual energy form. \vspace{1mm}

In this paper, we will focus on strongly local Dirichlet spaces. The purpose of this paper is to address the following questions: 
\begin{itemize}
	\item How can one define $\lambda$-subharmonicity (resp. $\lambda$-superharmonicity) for $L^1_\loc$--functions on such spaces, in a way that the definition is consistent with the distrubutional one in the Riemannian case? Note that this will lead naturally to the notion of an $L^q$-positivity preserving strongly local Dirichlet space. 
	\item What is the regularity of $\lambda$-subharmonic $L^1_\loc$--functions on strongly local Dirichlet spaces?
	\item Which class of strongly local Dirichlet spaces is $L^q$-positivity preserving?   
\end{itemize}

Addressing the first question, note first that on a general strongly local Dirichlet space there is no natural substitute for test functions. To deal with this problem, for every open $U\subset X$ let $P_t^U$ denote the semigroup of $(U,\m,\mathscr{E}, \F(U))$ with Dirichlet type boundary conditions. We define $f\in L^1_\loc(X,\m)$ to be \emph{locally $\lambda$-shift defective}, if for all open relatively compact sets $U\subset X$ there exists a weakly $\lambda$-harmonic function $g\in C_b(U)$ on $U$ such that 
$$
e^{-t\lambda}P_t^U(f|_U-g)\geq f|_U-g\quad\text{for all $t>0$.}
$$
This new notion, which can be considered as a localized and unsigned variant of the classical concept of defective functions, turns out on Riemannian manifolds to be equivalent to $\lambda$-subharmonicity in the distributional sense (cf. Proposition \ref{aspspo}).\vspace{1mm}

To proceed further with the theory, we introduce a new class of strongly local Dirichlet spaces, which we call \emph{locally smoothing spaces}. This assumption is entirely local, in the sense that all Riemannian manifolds (without any curvature assumptions) are locally smoothing and all open subsets of locally smoothing spaces again have this property. In addition, we show that all finite dimensional $\mathrm{RCD}$--spaces, all Carnot groups and the unbounded Sierpinski gasket belongs to this category of spaces. \vspace{1mm}

In order to establish a regularity theory for locally $\lambda$--shift defective functions, we assume throughout that $(X,\m,\mathscr{E}, \F(X))$ is locally smoothing. A central technical step, more or less intrinsic in the definition, is an approximation result, Lemma \ref{approx}, which reveals that locally $\lambda$--shift defective functions always have upper semicontinuous representatives. Furthermore, Theorem \ref{esqq} states that every locally $\lambda$--shift defective function in $\F_\loc(X)$ is weakly $\lambda$-subharmonic, and conversely, that every weakly $\lambda$-subharmonic which is locally bounded above is locally $\lambda$--shift defective (noting that on many spaces, every weakly $\lambda$-subharmonic is automatically locally bounded above; cf. Remark \ref{bemi}). Our main regularity result, Theorem \ref{main}, states that for every locally $\lambda$--shift defective function $f$ one has $(f-c)_+^{q/2}\in \F_\loc(X)$ for all $c\geq 0$, $q\in (1,\infty)$, which generalizes earlier results from \cite{sturm1}. The proof of this result is based on a Kato-Brezis type result, Theorem \ref{brezis}, which states that the pointwise maximum of two locally $\lambda$-shift defective functions has again this property. The proof of our Kato-Brezis theorem in turn relies on a probabilistic characterization of locally $\lambda$--shift defective functions, Proposition \ref{prob}, which is based on a probabilistic characterization of locally bounded weakly $\lambda$-subharmonic functions from \cite{kuwae}. We also obtain a strong maximum principle, Theorem \ref{th:maximumprinciple}, which does not need an $\F_\loc(X)$ or a continuity assumption. We note in passing that these results entail that if $f$ is a weakly $\lambda$-harmonic function on an $\mathrm{RCD}^*(K,N)$ space with $\lambda f\in L^q_\loc(X,\m)$ for some $q\geq \max(N,4)$, then the square $|\nabla f|^2_*$ of the minimal relaxed slope of $f$ admits on every open relatively compact subset an upper semicontinuous representative which is subject to the strong maximum principle (cf. Example \ref{rcdex}). Note that such a result is rather delicate, as on Alexandrov spaces it may happen that $|\nabla f|^2_*$ need not be continuous \cite{peng}, which follows essentially from the results obtained in \cite{filip,otso}.\vspace{1mm}

The regularity of harmonic and subharmonic functions is a central theme in analysis and beyond, ever since Weyl's seminal work \cite{weyl}. We refer to \cite{stroock} for a discussion of the latter paper and its impact, and to the textbooks \cite{evans,jost} for a general PDE point of view. It is also a classic fact that harmonicity properties come in semigroup resp. probabilistic flavors, see Kakutani's and Doob's fundamental papers \cite{kakutani}, \cite{doob} for the Euclidean space. It is no surprise that this aspect is very well captured within the framework of Dirichlet spaces, as those are in a 1-1 correspondence with Hunt processes \cite{fukushima_classic,maroeckner,fukushima,chenfukushima}. An important point of the intrigue is that different notions of (sub-)harmonicity require different a-priori regularity assumptions, while the theory from this paper is based on only an $L^1_\loc$--assumption. In this context, we note that recently a distributional notion of $\lambda$-subharmonic functions for finite dimensional RCD--spaces has been given in \cite{peng} (see in particular also \cite{gmss} for very interesting equivalent characterizations of this notion under additional local regularity assumptions). The theory from \cite{peng, gmss} uses the fact that a natural space of test functions is available on RCD-spaces. However, when applied to a complete Riemannian manifold with Ricci curvature bounded from below, this class of test functions is larger than smooth compactly supported functions, so that our definition of the distributional inequality $\Delta f\leq \lambda  f$ seems to be slightly more general even in this setting. \vspace{1mm}

Finally, based on our notion of locally $\lambda$-shift excessive functions, we define $L^q$-positivity preserving strongly local Dirichlet spaces, and prove the $L^q$--variant, $q\in (1,\infty)$, of the BMS conjecture for irreducible and intrinsically complete locally smoothing spaces. Our proof combines the $L^q$-Liouville property from \cite{sturm1} with our main regularity result, and in particular, does not rely on Laplacian cut-off functions, which neither make sense nor are available in this generality (note, however, that Laplacian cut-off functions exist in principle on finite dimensional RCD-spaces; cf. Lemma 3.1 in \cite{mondino}). We remark also that recently it has been shown in \cite{piggio} that in the Riemannian case, the $L^q$--variant of the BMS conjecture remains true, if one removes from a complete Riemannian manifold a possibly singular set with Hausdorff co-dimension strictly larger than $2q/(q-1)$, the threshold value $2q/(q-1)$ being sharp. \vspace{1mm}

The present paper is organized as follows: in Section \ref{sec1} standard notions and notation concerning local Dirichlet spaces, their associated operators and heat semigroups and the corresponding diffusion processes is recorded. In Section \ref{sec2} we introduce locally smoothing spaces and list a number of classes of strongly local Dirichlet spaces that fall into this framework. In Section \ref{sec3} we first review some important auxiliary results on weakly $\lambda$-subharmonic and $\lambda$-defective (signed) functions. We then give an equivalent characterization of distributionally $\lambda$-subharmonic functions on Riemannian manifolds in terms of the local heat semigroups, which motivates our notion of locally $\lambda$--shift defective functions, introduced next. We go on to prove that on locally smoothing spaces, the concepts \lq weak $\lambda$-subharmonicity' and \lq local $\lambda$-shift defectivity' are consistent under natural a priori regularity assumptions. Section \ref{sec4} contains the aforementioned probabilistic characterization of locally $\lambda$-shift defective functions, the Kato-Brezis theorem, and the regularity Theorem \ref{main}. Section \ref{sec5} is devoted to the maximum principle. Finally, Section \ref{sec6} contains a formulation and proof of the BMS-conjecture in our general framework.\vspace{3mm}

\emph{Acknowledgements:} The authors would like to thank Sebastian Boldt, Kazuhiro Kuwae and Marcel Schmidt for very helpful discussions. We are indebted to the referee for his careful reading and suggestions that improved the article.\vspace{3mm}

\emph{Conflicts of interest statement:} In this research there is no conflicts of interests.\vspace{3mm}

\emph{Data availability statement:} In this research no data is processed.

\section{Strongly local Dirichlet spaces}\label{sec1}

In the sequel we understand our function spaces over $\IR$. We follow standard notation and refer to \cite{fukushima,maroeckner} for textbooks on Dirichlet forms and \cite{dello1} for a recent article with an in depth treatment of Dirichlet forms.
  
\begin{definition} A \emph{strongly local Dirichlet space} is a quadruple $(X,\m,\mathscr{E}, \F(X))$ where $X$ is a locally compact, separable, metrizable space, equipped with a Radon measure $\m$ having full support, and $\mathscr{E}$ is a closed, symmetric, nonnegative bilinear form in $ L^2(X,\m)$ with dense domain of definition $\F(X)\subset L^2(X,\m)$, such that
\begin{itemize}
\item (Markovian property) for all contractions $T:\IR\to\IR$ with $T(0)=0$ and all $f\in\F(X)$ one has $T\circ f\in\F(X)$ with 
$$
\En(T\circ f,T\circ f)\leq \En(f,f),
$$ 
\item (Regularity) The space $\F(X)\cap C_c(X)$ is dense in $C_c(X)$ with respect to $\left\|\bullet\right\|_{\infty}$ and in $\F(X)$ with respect to the scalar product 
$$
\En_1(u,v):=\En(u,v)+\left\langle u,v \right\rangle,
$$
where
$$
\left\langle u,v\right\rangle=\int_{X} u v \Id\m,
$$
the scalar product in $L^2(X,\m)$, a notation we will use, more generally, as long as $uv\in L^1(X,\m)$.
\item (Strong locality) For all $u, v\in\F(X)$ such that $u$ is constant on the support of $v$:
$$\En(u,v)=0 .$$ 
\end{itemize}
\end{definition}

We fix a strongly local Dirichlet space $(X,\m, \mathscr{E}, \mathscr{F}(X))$.\vspace{1mm}

For open subsets $U,V\subset X$ we use the notation $U\Subset V$ to indicate that $\overline{U}$ is a compact subset of $X$ with $\overline{U}\subset V$.\vspace{1mm}

The regularity property implies \cite{fukushima} that $\mathscr{F}(X)\cap L^{\infty}(X,\m)$ is an algebra, and that for each $U_1\Subset U_2\Subset X$ there exists $\psi\in \mathscr{F}(X)\cap C_c(X)$ with $\psi=1$ on $\overline{U_1}$ and $\mathrm{supp}(\psi)\subset U_2$. From regularity it also follows that there is a Choquet capacity \cite{fukushima} $\mathrm{Cap}(\cdot)$ on $X$ such that for all open $U\subset X$ one has
$$
\mathrm{Cap}(U)=\inf\{\En_1(f)\;  | \; f\in\F(X), 1_U\leq f \}.
$$
One has $\m\ll \mathrm{Cap}$ on Borel sets, and properties that hold away from a set of zero capacity are said to hold quasi everywhere (q.e). A function $f:X\to\IR$ is called \emph{quasi-continuous}, if for all $\epsilon>0$ there exists an open set $U_\epsilon\subset X$ with $\mathrm{Cap}(U_\epsilon)\leq \epsilon$ such that $f|_{X\setminus U_\epsilon}$ is continuous. Every element of $\F(X)$ has a quasi-continuous representative and two such representatives agree q.e., \cite{fukushima}. 
Following \cite{sturm1}, see \cite{fukushima} as well, let
\begin{align*}
\efloc(X)&:=\set{f\in L^2_\loc(X,\m)}{\text{for all $U \Subset X$ there exists $h\in\F(X)$ with $f=h$ on $U$}}\\
& =\set{f\in L^2_\loc(X,\m)}{\text{ for all $\phi\in\F(X)\cap C_c(X)$ one has $\phi f\in \F(X)$}},
\end{align*}
with the above equality easily checked by regularity. By what we mentioned earlier, one has $1\in \mathscr{F}_\loc(X)$. Note that there are different ways of introducing $\efloc(X)$, see \cite{kuwae-98,frank} and the discussion in \cite{dello1}. We also note that thanks to locality, $\En$ has a canonical extension to a bilinear map $\efloc(X)\times \F_c(X)\to\mathbb{R}$, where $\F_c(X)$ denotes the elements of $\F$ with compact $\m$-support, given by $\En(f,\phi):=\En(h,\phi)$, where $h\in \F(X)$ is such that $h=f$ in an open neighbourhood of the $\m$-support of $\phi$.\vspace{1mm}

Every element of $\F_\loc(X)$ has a quasi-uniquely determined quasi-continuous representative, which we will always denote by the same symbol again.\vspace{1mm}

Strong locality implies the formula of \emph{Beurling-Deny and Le Jan}
$$
 \En( f,g)=\int_X \Id\Gamma(f,g)\quad\text{for all }f,g\in\F(X),
$$
where the symmetric, nonnegative, bilinear map
$$
\Gamma:\mathscr{F}(X)\times \mathscr{F}(X)\longrightarrow \{\text{signed finite Radon measures on $X$}\},
$$
the so called \emph{energy measure of $(X,\m, \mathscr{E}, \mathscr{F}(X))$}, is defined through polarization by 
$$
\int_X \phi \ \Id\Gamma(f,f)= \En( f,\phi f)-\frac{1}{2}\En(f^2,\phi)\quad\text{for all  }f\in\F(X)\cap L^{\infty}(X,\m), \phi\in \F(X)\cap C_c(X).
$$
Whenever it makes sense, we are going to use the standard notation $\En(f):=\En(f,f)$ and $\Gamma(f):=\Gamma(f,f)$. The energy measure has the following properties:
\begin{itemize}
\item $\Gamma$ \textrm{does not charge sets of zero capacity} \cite{sturm1, fukushima},
\item $\Gamma$ is strongly local, in the sense that for all $U\subset X$ open, all $u\in \F(X)$ such that $u$ is constant on $U$, and all $v\in\F(X)$, one has $1_U \Id \Gamma(u,v)=0$, cf. \cite{sturm1},
\item $\Gamma$ satisfies the \textrm{Leibniz rule} \cite{sturm1}, in the sense that for all 
\begin{align*}
\Id\Gamma(uv,w)=u\Id\Gamma(v,w)+v\Id\Gamma(u,w)\quad\text{for all $u,v \in\F(X)\cap L^{\infty}(X,\m)$, $w\in\F(X)$,}
\end{align*}
\item $\Gamma$ satisfies the \textrm{chain rule} \cite{sturm1}, in the sense that 
\begin{align*}
\Id\Gamma(\eta\circ u,v)=(\eta'\circ u)\Id\Gamma(u,v)\quad\text{for all $u,v\in\F(X)$,}
\end{align*}
where $\eta$ is any function which is $C^1$ with a bounded derivative, on an interval which contains the $\m$-essential image of $u$.
\item $\Gamma$ satisfies the \textrm{truncation property} \cite{sturm1}, in the sense that 
\begin{align*}
\Id\Gamma(u\wedge v,w)=1_{\{ u \geq v\}}\Id\Gamma(u,w)+1_{\{ u<v\}}\Id\Gamma(v,w)\quad\text{for all }u,v,w\in\F(X)
\end{align*}
\item $\Gamma$ satisfies the \textrm{Cauchy-Schwarz inequality}
\begin{align}\label{CS}
	\int_X fg\Id\Gamma(v,w)\leq \sqrt{\int_X f^2\Id\Gamma(v)}\sqrt{\int_X g^2\Id\Gamma(w)} \quad\text{for all $v,w\in\F(X)$, $f,g\in \F(X)\cap L^\infty(X,\m)$.}
\end{align}
\end{itemize}

Again, using the locality of $\Gamma$, the energy measure extends to a nonnegative definite symmetric bilinear form
$$
\Gamma:\mathscr{F}_\loc(X)\times \mathscr{F}_\loc(X)\longrightarrow \{\text{signed Radon measures on $X$}\},
$$
namely, given $f_1,f_2\in\F_\loc(X)$ and $U\Subset X$ one sets
$$
1_U \Id \Gamma(f_1,f_2):= 1_U \Id \Gamma(h_1,h_2),
$$
where $h_1,h_2\in\F(X)$ are chosen such that $f_j=h_j$ on $U$. All the above listed properties of the energy measure remain to hold, if one replaces $\F(X)$ with $\F_\loc(X)$ and $L^{\infty}(X,\m)$ with $L^{\infty}_\loc(X,\m)$ accordingly in the statements. \vspace{2mm}

For future reference, we also note:

\begin{lemma}\label{peter} Let $0\leq w\in L^\infty_\loc(X,\m)$ and assume there exists a sequence $\seqx{w}{k}$ in  $\F_\loc(X)$ such that $w_k\to w$ in $L^2_\loc(X,\m)$ and
		$$
		\sup_{k}\int_X \mathrm{d}\Gamma(w_k)<\infty.
		$$
Then $w\in\F_\loc(X)$.
\end{lemma}

\begin{proof} It follows easily from the Markovian property that $0\leq v\in L^\infty_\loc(X,\m)$ is in $\F_\loc(X)$, if and only if for all $\phi\in C_c(X)\cap \F(X)$ one has $\phi\wedge v \in \F(X)$. Thus, fixing an arbitrary $\phi\in C_c(X)\cap \F(X)$, we have $\phi\wedge w_k\in\F(X)$ for all $k$, and it remains to show $\phi\wedge w \in \F(X)$. To this end, we note  $ \phi\wedge w_k\to \phi\wedge w$ in $L^2(X,\m)$. Moreover, by the truncation property,
$$
\En(w_k\wedge \phi)=\int_{\{\phi<w_k\}}\mathrm{d}\Gamma(\phi)+\int_{\{\phi\geq w_k\}}\mathrm{d}\Gamma(w_k) \leq \En(\phi)+\sup_k \int_X\mathrm{d}\Gamma(w_k).  
$$
It follows that $\seqx{\phi\wedge w}{k}$ is bounded in the Hilbert space $(\F(X),\En_1)$. Weak compactness implies that it has a weakly convergent subsequence. The limit must agree with $\phi\wedge w$, giving $\phi\wedge w \in \F(X)$.
\end{proof}	

We use the notation 
$$
\En_\lambda(u,v):=\En(u,v)+\lambda\left\langle u,v \right\rangle,
$$
whenever it makes sense, and recall that a function $f\in \mathscr{F}_\loc(X)$ is called 
\begin{itemize}
	\item \emph{weakly $\lambda$--subharmonic}, if 
\begin{align}\label{ert}
	\text{$\mathscr{E}_\lambda(f,\phi)\leq 0$ for all $0\leq \phi\in \mathscr{F}_c(X)$}.
\end{align}
\item \emph{weakly $\lambda$--superharmonic}, if $-f$ is weakly $\lambda$--subharmonic, 

\item \emph{weakly $\lambda$--harmonic}, if $f$ is weakly $\lambda$--subharmonic and weakly $\lambda$--superharmonic.
 \end{itemize}
If $\lambda=0$, one simply talks about weak (sub-/super-) harmonicity.

\begin{lemma}\label{sswwq} Assume $\lambda\geq 0$ and that $f_1,f_2\in \F_\loc(X)$ are weakly $\lambda$-subharmonic. Then $f_1\vee f_2$ is weakly $\lambda$-subharmonic.
\end{lemma}

\begin{proof} In view of $f_1\vee f_2\in \F_\loc(X)$, this statement follows from Theorem 6.4 (iii) in \cite{kuwaeMax}.
\end{proof}

Let $H\geq 0$ denote the nonnegative self-adjoint operator in $L^2(X,\m)$ induced by $\En$, \cite{kato,fukushima}, and let $P_t:=\mathrm{e}^{-tH}$ for $t\ge 0$ (defined by the functional calculus). We call $P=(P_t;t\ge 0)$ the associated (self-adjoint) heat semigroup in $L^2(X,\m)$. The heat semigroup extends to a positivity preserving Markovian contraction semigroup in $L^q(X)$ for all $q\in [1,\infty]$, which is consistent in $q$ (therefore we use the same symbol for the formally different semigroups on the different spaces), strongly continuous for $q<\infty$ and weak-$*$-continuous for $q=\infty$. \\
Let $\mathscr{B}(X)$ and $\mathscr{B}(X,\m)$ denote, respectively, the space of Borel functions on $X$ and the space of $\m$-equivalences of Borel functions on $X$. We stretch the notation further by setting 
$$
P_tf(x):=\lim_{n\to\infty} P_t(f\vee (-n))(x)\in [-\infty,0]\quad \text{for every $0\geq f\in \mathscr{B}(X,\m)$,}
$$
and
$$
P_tf(x):=\lim_{n\to\infty} P_t(f\wedge n)(x)\in [0,\infty]\quad \text{for every $0\leq f\in \mathscr{B}(X,\m)$.}
$$

Let $\infty_X\notin X$ be a cemetery point. We equip the space $\IWW(X)$ of continuous paths 
$$
[0,\infty)\longrightarrow  \hat{X}:=X\cup\{\infty_X\}
$$
taking values in the Alexandrov compactification $\hat{X}$ of $X$ with the $\sigma$-algebra and the filtration which is induced by the coordinate process
$$
\mathbb{X}: [0,\infty)\times \IWW(X)\longrightarrow \hat{X},
$$
with the usual notation $\mathbb{X}_s:=\mathbb{X}(s,\bullet)$. Let
$$
\theta_s:\IWW(X)\longrightarrow \IWW(X),\quad \gamma\longmapsto \gamma(s+\bullet)
$$ 
denote the shift operator, and, given an open subset $V\subset X$, let
$$
\tau_V:\IWW(X)\longrightarrow [0,\infty], \quad \tau_V(\gamma):=\inf\{s\geq 0: \mathbb{X}_s(\gamma)\in X\setminus V\}.
$$
It follows from a by now classical result of Fukushima \cite{fukushima_classic} that there exists a family $(\mathbb{P}^x)_{x\in X}$ of diffusion measures on $\IWW(X)$ such that for every $f\in \mathscr{B}(X)\cap L^2(X,\m)$ one has
\begin{equation}\label{process}
P_tf(x)=\mathbb{E}^x[f(\mathbb{X}_t)]:=\int f(\mathbb{X}_t) \mathrm{d}\mathbb{P}^x\quad\text{for all $t>0$, $\m$-a.e. $x \in X$},
\end{equation}
where $f(\infty_X):=0$. See also \cite{maroeckner} for the quasi-regular case. \\

Let $U\subset X$ be an open subset. Whenever convenient, we identify functions (and $\m-$equivalence classes) that vanish outside $U$ with their restriction to $U$. We get the strongly local Dirichlet space $(U,\m,\mathscr{E},\mathscr{F}(U))$ by setting
$$
\mathscr{F}(U)=\overline{C_c(U)\cap \F(X)}^{\En_1}=\set{f\in  \mathscr{F}(X)}{f=0 \mbox{  q.e in }X\setminus U}.
$$ 
The induced operator, semigroup and diffusion measures are denoted by $H^{U}$, $P^{U}$, $(\mathbb{P}_U^x)_{x\in X}$, respectively. It is well-known that for all open $U\subset V\subset X$ one has 
$$
P^U_tf \leq P^V_tf\quad\text{for all $t>0$, $0\leq f\in \mathscr{B}(X,\m)$}.
$$
Moreover, from monotone convergence theorems for forms, \cite{kato,simon}, it easily follows that 
$$
P^{U_n}_tf\uparrow P^V_tf\quad\text{as $n\to\infty$, for all $t>0$, $0\leq f\in \mathscr{B}(X,\m)$,}
$$
whenever $(U_n)_{n\in\IN}$ is a sequence of open subsets of $V$ with $U_n\uparrow V$.

\section{Locally smoothing spaces}\label{sec2}

We collect some classical notions, slightly adjusted to our situation, in:

\begin{definition} We say that [the semigroup of] $(X,\m,\En,\F(X))$ is
\begin{itemize}
\item \emph{ultracontractive}, if 
\begin{align}\label{ultra}
	P_t:L^1(X,\m)\longrightarrow L^\infty(X,\m)\quad\text{for all }t>0,
\end{align}
by which we mean that $P_t(L^1(X,\m))\subset L^\infty(X,\m)$ (which in turn gives continuity of the induced operator by the closed graph theorem),
\item \emph{doubly Feller}, if one has the mapping properties
\begin{align}\label{feller}
P_t:C_0(X)\longrightarrow C_0(X),\quad \text{for all }t\geq 0,
\end{align}
and
$$
P_t:L^\infty(X,\m)\longrightarrow C(X)\quad\text{for all }t>0,
$$
and in addition
\begin{align}\label{pointwise}
\lim_{t\to 0+}P_tf(x)=f(x)\quad\text{for all $f\in C_0(X)$, $x\in X$.}
\end{align}
Here  $C_0(X)$ denotes the family of
continuous functions vanishing at infinity.
\item \emph{irreducible}, if every $f\in\mathscr{F}_\loc(X)$ with $\mathscr{E}(f)=0$ (which is equivalent to $\Gamma(f)=0$) is constant.

\end{itemize}
 We say that [the generator of] $(X,\m,\En,\F(X))$ 
\begin{itemize}
\item \emph{has a spectral gap}, if $H>0$, i.e., $\inf(\sigma(H))>0$.\vspace{1mm}

\item \emph{has purely discrete spectrum}, if the spectrum $\sigma(H)$ of $H$ is purely discrete, that is, consists of isolated eigenvalues having a finite multiplicity. 
\end{itemize}
\end{definition}

A few simple remarks are in order:

\begin{remark}\label{locbem} 1. Assume (\ref{feller}) and (\ref{pointwise}). 
Then, by Exercise I-(9.13), p. 51 in \cite{blumen} or Lemma 1.4 in \cite{boetcher}, one gets that $P_t:C_0(X)\to C_0(X)$, $t>0$, is strongly continuous. In particular, then one also has the mapping property
$$
(H+\lambda)^{-1}:C_0(X)\longrightarrow C_0(X),\quad \text{for all }\lambda>0,
$$
which follows from the representation
\begin{align}\label{laplace}
(H+\lambda)^{-1}f=\int^\infty_0 e^{-\lambda t}P_tf dt,\quad f\in C_0(X),
\end{align}
where the improper Riemannian integral converges in the uniform norm.\vspace{1mm}

2. If $(X,\m,\En,\F(X))$ has purely discrete spectrum, then $(X,\m,\En,\F(X))$ has a spectral gap, if and only if $\mathrm{Ker}(H)=\{0\}$.\vspace{1mm}

3. The space $(X,\m,\En,\F(X))$ is irreducible \cite{sturm1}, if and only if for every Borel set $Y\subset X$ with $P_t(1_Yf)=1_YP_tf$ for all $t>0$, $f\in L^2(X,\m)$ one has either $\m(Y)=0$ or $\m(X\setminus Y)=0$. This is equivalent to the following positivity improving property: for all $t>0$, $0\leq f\in L^2(X,\m)\setminus \{0\}$ one has $P_tf>0$. \vspace{1mm} 

4. If $(X,\m,\En,\F(X))$ is irreducible, then for all $0\leq f\in L^2(X,\m)\setminus\{0\}$, $\lambda>0$ one has $(H+\lambda)^{-1}f>0$. This follows from the previous remark and formula (\ref{laplace}), which converges for $f \in L^2(X,\m)$ in the $L^2$--Bochner sense.\vspace{1mm}

5. By self-adjointness, (\ref{ultra}) is equivalent to 
\begin{align*}
	P_t:L^1(X,\m)\longrightarrow L^2(X,\m)\quad\text{for all $t>0$,}
\end{align*}
in particular, by the spectral theorem and the semigroup property one then has
$$
P_t:L^1(X,\m)\longrightarrow \bigcap_{n\in\IN}\dom(H^n)\quad\text{for all }t>0.
$$
\end{remark}

We record that in the doubly Feller case, the underlying diffusion is pointwise determined by the semigroup:

\begin{lemma} Assume that $(X,\m,\En,\F(X))$ doubly Feller. Then there exists a uniquely determined family of diffusion measures $(\mathbb{P}^x)_{x\in X}$ on $\IWW(X)$, such that for every bounded function $f:X\to\IR$, every $t>0$, and \emph{every} $x\in X$, one has
\begin{align}\label{formi}
P_tf(x)=\mathbb{E}^x[f(\mathbb{X}_t)].
\end{align}
Moreover, one has the absolute continuity condition 
$$
(\mathbb{X}_t)_*\mathbb{P}^x\ll\m\quad\text{for all $t>0$, $x\in X$,}
$$
in particular, for all $q\in [1,\infty]$, $f\in L^q(X,\m)$, $t>0$, the function $P_tf$ has a pointwise well-defined $\m-$representative, given by the RHS of (\ref{formi}).
\end{lemma}

\begin{proof} 
This result is standard, except possibly the fact that $\mathbb{P}^x$ is concentrated on continuous paths for \emph{every} $x\in X$ (noting that the standard theory of Feller processes only produces measures that are concentrated on cadlag paths). This above pointwise concentration follows from observing that the strong locality of $\mathscr{E}$ implies this property for $\m$-a.e. $x$ (in fact, q.e. $x\in X$), and then the doubly Feller property together with Theorem 4.5.4(iii) in \cite{fukushima} allows to conclude the pointwise result.
\end{proof}

\begin{definition}\label{raum} $(X,\m,\mathscr{E}, \mathscr{F}(X))$ is called \emph{locally smoothing}, if for all $U\Subset X$ it holds that
\begin{itemize}
\item $(U,\m,\En,\F(U))$ is ultracontractive,
\item there exists $U'\subset X$ open with $U\subset U'$ and $U'\setminus \overline{U}\neq \emptyset$ such that $(U',\m,\En,\F(U'))$ is doubly Feller and irreducible.   	
\end{itemize}	
\end{definition}

Note that if $(X,\m,\mathscr{E}, \mathscr{F}(X))$ is locally smoothing, then so is $(U,\m,\mathscr{E}, \mathscr{F}(U))$ for every open $U\subset X$.\vspace{1mm}

Locally smoothing spaces have the following self-improvement properties:

\begin{proposition}\label{oooo} Assume that $(X,\m,\mathscr{E}, \mathscr{F}(X))$ is locally smoothing, and let $U\Subset X$.\\
\emph{a)} For all $t>0$ one has 
\begin{align}\label{asdfg}
	P^{U}_t:L^1(U,\m)\longrightarrow C_b(U).
\end{align} 
\emph{b)} The generator of $(U,\m,\En,\F(U))$ has purely discrete spectrum and a spectral gap.
\end{proposition}

\begin{proof}(a) This statement is well-known in various forms (cf. \cite{ck,chung2}). We give a detailed proof for the convenience of the reader: by the semigroup property, we only have to show that $P^U_tf$ is continuous for all $f\in L^\infty(U,\m)$. To see this, pick $U\subset U'$ open with $(U',\m,\En,\F(U'))$ doubly Feller. We first record the standard fact that for all $t>0$, $\m$-a.e. $x\in U$ one has 
	$$
	P^{U}_tf(x)= \mathbb{E}^x_{U'}[1_{\{t<\tau_U\}}f(\mathbb{X}_t)],
	$$
	and we consider the LHS to be defined pointwise by the RHS. It follows that for all $0<s<t$ one has, by a standard calculation that relies on the Markov property (cf. p. 18 in \cite{schilling}), 
	$$
	P^{U}_tf(x)-\big(P^{U'}_sP^{U}_{t-s}f \big)(x)= \mathbb{E}^x_{U'}\left[1_{\{s\geq \tau_U\}}\theta_s(1_{\{t-s<\tau_U\}}f(\mathbb{X}_{t-s}))\right]\quad\text{for all }x\in U.
	$$ 
Thus,
	$$
	\left|P^{U}_tf(x)-\big(P^{U'}_sP^{U}_{t-s}f \big)(x)\right|\leq \left\|f\right\|_{\infty}\IP^x_{U'}(s\geq \tau_U).
	$$
	Since $P^{U'}_sP^{U}_{t-s}f$ is continuous, it remains to show that $\IP^x_{U'}(s\geq \tau_U)$ converges locally uniformly in $x$ to $0$ as $s\to 0$. This, however, is a well-known consequence of the Feller property (\ref{feller}) of $(U',\m,\En,\F(U'))$ (cf. Lemma 2.5 in \cite{chung}).\vspace{1mm}
	
b) Pick $U\subset U'$ open with  $U'\setminus \overline{U}\neq \emptyset$ such that $(U',\m,\En,\F(U'))$ is irreducible. Since $P^U_t$ maps to $L^\infty(U,\m)$ and $\m(U)<\infty$, $P^U_t$ is Hilbert-Schmidt for every $t>0$, see \cite{df}, 11.2, 11.16 and the discussion in \cite{stollman}, p.418; thus $H^U$ has purely discrete spectrum. Assume $\varphi_U\in \mathrm{Ker}(H^U)$. We are going to show $\varphi_U=0$. As we have $\varphi_U\in \mathscr{F}(U)\subset\mathscr{F}(U')$ with $\En(\varphi_U)=0$, it follows from the irreducibility of $(U',\m,\En,\F(U'))$ that $\varphi_U$ is constant $\m$-a.e. on $U'$. On the other hand, $\m(U'\setminus U)>0$ and we have $\varphi_U=0$ $\m$-a.e. in $U'\setminus U$, thus $\varphi_U=0$.

\end{proof}

Let us give some examples of spaces that are locally smoothing:

\begin{example}\label{esxcv} \textbf{Arbitrary Riemannian manifolds:} Let $X$ be a connected Riemannian manifold with its Levi-Civita connection $\nabla$ and its Riemannian volume measure $\m$. Then with
\begin{align}\label{syp}
\mathscr{E}(f_1,f_2):=\int_X (\nabla f_1, \nabla f_2) \mathrm{d}\m 
\end{align}
the usual energy form with domain of definition $\mathscr{F}(X)=W^{1,2}_0(X)$, the triple $(X,\m,\En,\F(X))$ is well-known to be a strongly local Dirichlet space. In this case one has $\mathscr{F}_\loc(X)= W^{1,2}_\loc(X)$, and $\mathscr{F}(U)= W^{1,2}_0(U)$ for all open $U\subset X$, so that the restricted strongly local Dirichlet spaces correspond to Dirichlet boundary conditions. Moreover, $(X,\m,\En,\F(X))$ is locally smoothing: indeed, given $U \Subset X$, one has (\ref{ultra}) since for the heat kernels one has
$$
p^U(t,x,y)\leq p(t,x,y),
$$ 
and $p(t,x,y)$ is jointly smooth in $(x,y)$ (in fact in $(t,x,y)$). Moreover, picking $U'\Subset X$ connected and with smooth boundary such that $U\subset U'$ and $U'\setminus \overline{U}\neq \emptyset$ one has that $(U',\m,\En,\F(U'))$ is doubly Feller, because $p^{U'}(t,x,y)$ is smooth, and $P_t^{U'}f$ extends continously to zero on $\partial U'$ because of Dirichlet boundary conditions (the strong continuity of $P^{U'}_t:C_0(U')\to C_0(U')$ is automatic in this case). Moreover, $(U',\m,\En,\F(U'))$ is irreducible, because the connectedness of $U'$ implies that $p^{U'}(t,x,y)>0$.
\end{example}

\begin{example}\label{rcd} \textbf{RCD$^*$--spaces:} Consider a complete, locally compact, separable metric measure space $(X,\dm)$, and a Radon measure $\m$ on $X$ with full support that gives a finite mass to open balls; a function $g\in L^2(X,\m)$ is called a \emph{relaxed slope of $f\in L^2(X,\m)$}, if there exists $\tilde{g}\in  L^2(X,\m)$ and a sequence of Lipschitz functions $f_n\in L^2(X,\m)$ such that 
\begin{itemize}
	\item $f_n\to f$ in $L^2(X,\m)$ and $\mathrm{Lip} (f_n)\to \tilde{g}$ weakly in $L^2(X,\m)$
	\item $\tilde{g}\leq g$,
\end{itemize}
where
$$
\mathrm{Lip} (h)(x):=\limsup_{y\to x}\frac{|h(x)-h(y)|}{\dm(x,y)}
$$
denotes the local Lipschitz constant of a local Lipschitz function $h:X\to\IR$. Such a $g$ is called the \emph{minimal relaxed slope of $f$}, if its $L^2$-norm is minimal amongst all relaxed slopes of $f$, and then one sets $|\nabla f|_*:=g$. We refer the reader to \cite{ambro} for equivalent definitions of $|\nabla f|_*$. \\
With
$$
\F(X):=\{f\in L^2(X,\m):\text{$f$ has a relaxed slope}\},
$$
the \emph{Cheeger form} $\En$ is the densely defined functional on $L^2(X,\m)$ given by
$$  
\En(f):=\int_X |\nabla f|_*^2\mathrm{d}\m,\quad f\in \F(X),
$$
and $(X,\dm, \m)$ is called \emph{infinitesimally Hilbertian}, if $\En$ is a quadratic form. In this case, given $f,g\in \F(X)$ the limit
$$
(\nabla f,\nabla g):=\lim_{\epsilon\to 0+}\frac{1}{2\epsilon}\left(|\nabla (f+\epsilon g)|^2_*-|\nabla f|^2_*\right)
$$
exists in $L^1(X,\m)$, and $(X,\m,\En,\F(X))$ becomes a strongly local Dirichlet space, with 
$$
\mathrm{d}\Gamma(f,g)= (\nabla f,\nabla g)\mathrm{d} \m\quad\text{for all $f,g\in \F(X)$.}
$$
In the Hilbertian case, given $N\geq 1$, $K\in\IR$, one calls $(X,\dm,\m)$ an $\mathrm{RCD}^*(K,N)$-space (see \cite{eks} and the references therein), if 
\begin{itemize}
	\item every $f\in \F(X)$ with $|\nabla f|_*\leq 1$ has a $1$-Lipschitz $\m$-representative,
	\item one has
	$$
	\int_X e^{-\dm(x_0,x)^2}\mathrm{d}\m(x)<\infty\quad\text{for all $x_0\in X$, $c>0$,}
	$$
	\item for all $f\in\dom(H)$ with $Hf\in \F(X)$ and all $0\leq g \in \dom(H)\cap L^\infty(X,\m)$ with $Hg\in L^\infty(X,\m)$ one has the \emph{Bochner inequality}
$$
-\frac{1}{2}\int_X (H g)|\nabla f|^2_* \mathrm{d}\m+\int_X g\cdot (\nabla (Hf),\nabla f)\mathrm{d}\m\geq K\int_X g |\nabla f|^2_* \mathrm{d}\m+\frac{1}{N}\int_X g \cdot (H f)^2\mathrm{d}\m.
$$
\end{itemize}
Here $H$ denotes the nonnegative self-adjoint
operator in $L^2
(X, \m)$ associated with $(\En, \F(X))$.
If $(X,\dm,\m)$ is an $\mathrm{RCD}^*(K,N)$--space, then this space satisfies the local volume doubling assumption \cite{eks, sturm2}, and the heat kernel of $P_t$ satisfies the two-sided local Li-Yau heat kernel bounds \cite{jiang}. It follows from Proposition \ref{app} in the appendix that then $(X,\m,\En,\F(X))$ is locally smoothing.\\
Note that if $(X,\dm,\m)$ is induced from a complete connected Riemannian manifold $X$, then the Cheeger form is equal to the energy form, and the $\mathrm{RCD}^*(K,N)$-assumption is equivalent to $\dim(X)\leq N$ and $\mathrm{Ric}\geq K$.
\end{example}

\begin{example}\label{carnot} \textbf{Carnot groups:} Assume $G$ is a simply connected Lie group with Lie algebra $\mathfrak{g}$, such that there exists $N\geq 1$ and a stratification
	$$
	\mathfrak{g}=	\mathfrak{V}_1\oplus \cdots\oplus  \mathfrak{V}_N
	$$
	such that $[\mathfrak{V}_i,\mathfrak{V}_j]=\mathfrak{V}_{i+j}$ (noting that such stratification is essentially unique \cite{ledonne}), where $\mathfrak{V}_k:=\{0\}$ for $k>N$. In this situation, $G$ is nilpotent and called a \emph{Carnot group}. Let $\m$ denote the Haar measure on $G$, and let $V_1,\dots,V_d$ denote a basis of $\mathfrak{V}_1$, considered as left invariant vector fields on $G$. Then the closure $(\En,\F(G))$ in $L^2(G,\m)$ of the symmetric bilinear form
	$$
	C^\infty_c(G)\times C^\infty_c(G)\ni (f,g)\longmapsto \int_G \sum^d_{i,j=1}(V_i f) (V_j g) \mathrm{d} \m\in\IR
	$$ 
	turns $(X,\m,\En,\F(X))$ into a strongly local Dirichlet space with
	$$
	\mathrm{d}\Gamma(f,g)= \sum^d_{i,j=1}(V_i f) (V_j g)\mathrm{d} \m\quad\text{for all $f,g\in \F(G)$.}
	$$
	The Carnot-Caratheodory distance $\dm$ on $G$ is complete and induces the original topology on $G$ \cite{dragoni}, and the metric measure space $(G,\dm,\m)$ becomes doubling and Ahlfors regular, $\m(\B(x,r))\sim r^Q$, where
	$$
	Q:=\sum^N_{i=1}i \dim(\mathfrak{V}_i)
	$$
	denotes the homogeneous dimension. Furthermore, the heat kernel of $P_t$ satisfies the two-sided Gaussian estimate, as shown in \cite{bonf}: there exist $c_1,c_2,c_3,c_4>0$ such that for all $t>0$, $x,y\in G$ one has
	$$
	c_1 t^{-Q/2} e^{-c_2\frac{\dm(x,y)^2}{t}}\leq p(t,x,y)\leq 	c_3 t^{-Q/2} e^{-c_4\frac{\dm(x,y)^2}{t}}.
	$$
	It follows from Proposition \ref{app} in the appendix that $(G,\m,\En,\F(G))$ is locally smoothing. Note that in general $(G,\dm,\m)$ is not an $\mathrm{RCD}^*$-space. In fact, it has been shown in \cite{juillet} that for the Heisenberg group $\mathbb{H}_{2m+1}$ there exists no pair $(N,K)$ such that $(\mathbb{H}_{2m+1},\dm,\m)$ becomes an $\mathrm{RCD}^*(K,N)$ space.  
	\end{example}

\begin{example}\label{fractal} \textbf{Unbounded Sierpinski gasket:} Let $X\subset \IR^2$ be the unbounded Sierpinski gasket, with its usual graph approximation $(X_n)\subset X$. Then with $\alpha: =\log(3)/\log(2)$ the Hausdorff dimension of $X$ and $\m$ the $\alpha$-dimensional Hausdorff measure on $X$, the symmetric bilinear form
	$$
	\En(f,g):=\lim_{n\to\infty}\left(\frac{5}{3}\right)^n\sum_{x,y\in X_n, x\sim y} (f(x)-f(y))^2,\quad \F(X):=\{f\in L^2(X,\m): \En(f)<\infty\},
	$$
	turns $(X,\m,\En,\F(X))$ into a strongly local Dirichlet space. This follows from combining the results in \cite{barlow} with Proposition \ref{app2}.
\end{example}

\section{A new notion of subharmonic functions}\label{sec3}

We start this section with the following classical definition, see \cite{fukushima,maroeckner}:

\begin{definition} Given $\lambda\geq 0$, $h\in L^1_\loc(X,\m)$ is called 
	\begin{itemize}
		\item[(i)] \emph{$\lambda$--defective}, if $h\leq 0$ and $e^{-\lambda t}P_t h\geq h$ for all $t>0$,
		\item[(ii)] \emph{$\lambda$--excessive}, if $-h$ is $\lambda$--defective.
	\end{itemize}
\end{definition}

The following characterization of subharmonic functions of fixed sign is well-known, if one considers functions in $\F_\loc(X)$, see \cite{fukushima}, Theorem 2.2.1, \cite{maroeckner}, Proposition III.1.2. The main point of the result below is that we can drop this local regularity assumption. For $\lambda=0$, it has already been recorded by Sturm in \cite{sturm1}, Lemma 3, that this is possible. Note, however, that the methods from therein do not generalize directly to the $\lambda>0$ case, roughly speaking, as they rely on the strong locality of $\mathscr{E}$ (which fails for $\mathscr{E}_\lambda$).

\begin{theorem}\label{esqq0} Let $h\in L^1_\loc(X,\m)$, and let $\lambda\geq 0$.\\
\emph{(a)} Assume $h\leq 0$ and consider the following properties:\\
		\emph{(i)} $h$ is weakly $\lambda$--subharmonic.\\
		 \emph{(ii)} For all $f\in \F_c(X)$ one has $h\vee f\in\F(X)$ and 
		$\En_\lambda(h\vee f,f-h\vee f)\ge 0.$ \\
		 \emph{(iii)} For all $f\in \F_c(X)$ one has $h\vee f\in\F(X)$ and 
		$\mathscr{E}_\lambda(f,f-h\vee f)\ge 0.$ \\  
		\emph{(iv)}  For all $f\in \F_c(X)$ one has $h\vee f\in\F(X)$ and 
		$\mathscr{E}_\lambda(h\vee f)\le \mathscr{E}_\lambda(f).$\\
		\emph{(v)}  $h$ is  $\lambda$--defective.\\
	Then the following implications hold:  
		\emph{(i)}$\Longrightarrow$\emph{(ii)}$\Longrightarrow$\emph{[(iii)}\emph{ and (iv)]}$\Longleftrightarrow$\emph{(v)}.\\
\emph{(b)} Assume $h$ is $\lambda$--defective. Then one has
\begin{align*}
h\in\efloc(X)&\Longleftrightarrow\exists h^\sharp\in \efloc(X): h\ge h^\sharp\\
&\Longleftrightarrow\forall V\Subset X \ \exists h^\sharp_V\in\F_c(X): h|_V\ge h^\sharp_V.
\end{align*}
In particular, if $h\in L_\loc^\infty(X,\m)$, then $h\in\efloc(X)$.\\
\emph{(c)} If $h$ is $\lambda$--defective and $h$ satisfies one of the equivalent conditions from \emph{(b)}, then $h$ is weakly $\lambda$--subharmonic.
\end{theorem}

\begin{proof} We follow an approach which is in the spirit of \cite{schmidt}, Section 2.4.\\
(a) (i)$\Longrightarrow$(ii): This is \lq algebraic' in nature, we follow the proof of Lemma 2.50 in \cite{schmidt}. We let   $f\in \F_c(X)$ and consider $V\Subset X$ so that $\supp (f)\subset V$. Note that $h\in\efloc (X)$ by assumption and therefore, there is $h^\flat\in \F_c(X)$ such that 
	$$
	h=h^\flat\mbox{  on  }V .
	$$
	We can and will assume that $h^\flat\le 0$ and note that
	$$
	h\vee f=h^\flat \vee f \in \F_c (V) .
	$$
	Having seen this, we omit the superscript $^\flat$ for the rest of the calculation. 
	Since
	$$
	h\vee f =\frac12\left(h+f + \left| h-f\right|\right),
	$$
	it follows that
	\begin{align*}
		\mathscr{E}_\lambda(h\vee f)- \mathscr{E}_\lambda(h\vee f, f)
		&=\mathscr{E}_\lambda(h\vee f, h\vee f -f)\\
		&= \frac{1}{4}\mathscr{E}_\lambda\left( f+h + \left| h-f\right|, h  + \left| h-f\right|-f \right)\\
		&= \frac{1}{4}\left\{\mathscr{E}_\lambda\left(h + \left| h-f\right|\right)-\mathscr{E}_\lambda\left(f\right)\right\}\\
		&= \frac{1}{4}\left\{\mathscr{E}_\lambda\left(h\right)+ \mathscr{E}_\lambda\left(\left| h-f\right|\right) + 2\mathscr{E}_\lambda\left(h,\left| h-f\right|\right)-\mathscr{E}_\lambda\left(f \right)\right\}\\
		&\le \frac{1}{4}\left\{\mathscr{E}_\lambda\left(h\right)+ \mathscr{E}_\lambda\left( h-f\right) + 2\mathscr{E}_\lambda\left(h,\left| h-f\right|\right)-\mathscr{E}_\lambda\left(f \right)\right\}
		\\
		&= \frac{1}{4}\left\{2\mathscr{E}_\lambda\left(h\right)- 2\mathscr{E}_\lambda\left(h, f\right)+ 2\mathscr{E}_\lambda\left(h,\left| h-f\right|\right)\right\}\\
		&= \frac{1}{2}\left\{\mathscr{E}_\lambda\left(h,h-f + \left| h-f\right|\right)\right\}\\
		&= \mathscr{E}_\lambda\left(h,(h-f)_+\right)\le 0,
	\end{align*}
	since $h$ is weakly $\lambda$--subharmonic and
		$$
	0\le (h - f)_+=h\vee f - f\in \F_c(X) .
	$$
	From this chain of inequalities we get:
	\begin{equation*}
		\mathscr{E}_\lambda(h\vee f, f)\ge \mathscr{E}_\lambda(h\vee f) ,
	\end{equation*}
	which is (ii); expanding $f=f-h\vee f+h\vee f$ we get (iii), since 
	$$\mathscr{E}_\lambda(f,f-h\vee f)\ge\mathscr{E}_\lambda(h\vee f,f-h\vee f)\ge 0,$$
	which in turn is obvious. From (ii) and (iii) we get (iv):
 \begin{align*}
 \mathscr{E}_\lambda(f)-\mathscr{E}_\lambda(h\vee f)&= \mathscr{E}_\lambda(f+h\vee f, f-h\vee f)\\
 &=\mathscr{E}_\lambda(f,f-h\vee f)+
 \mathscr{E}_\lambda(f\vee h,f-h\vee f)\\
 &\ge 0
	\end{align*}
	The equivalence [(iii) and (iv)]$\Longleftrightarrow$(v) is a consequence of the characterization of convex sets that are invariant under the semigroup in terms of the corresponding quadratic form given in \cite{ouhabaz}. Here are the details: For given $h$ as in the assumption, the set
	$$
	K:=\{v\in L^2(X,\m)\; | \; v\ge h\} 
	$$
	is closed and convex in $L^2(X,\m)$ and the projection $\pi:L^2(X,\m)\to K$ is given by $\pi f=h\vee f$. 
Clearly (v) implies that the semigroup associated with $\mathscr{E}_\lambda$, given by $S_t=e^{-\lambda t}P_t$ for $t\ge 0$, leaves $K$ invariant.
Denoting the following stronger versions of the properties in the Theorem by\\
		 {(iii')} For all $f\in \F(X)$ one has $h\vee f\in\F(X)$ and 
		$\mathscr{E}_\lambda(f,f-h\vee f)\ge 0.$ \\  
		{(iv')}  For all $f\in \F(X)$ one has $h\vee f\in\F(X)$ and 
		$\mathscr{E}_\lambda(h\vee f)\le \mathscr{E}_\lambda(f),$\\
(iii')$\Longleftrightarrow$(iv')$\Longleftrightarrow$(v) follows from Theorem 2.1 and Corollary 2.4 in \cite{ouhabaz}.

 A standard approximation argument gives (iv) $\Rightarrow$ (iv') and so we get the implications: 
 
 (ii) $\Rightarrow$ [(iii) and (iv)] $\Rightarrow$ (iv') $\Leftrightarrow$ (v) $\Leftrightarrow$ (iii'), which settles (a).

 \vspace{1mm}
	
	(b) For the two equivalences, it remains to show that the last property implies the first. The other implications are evident. Let $ V\Subset X$ and $h^\sharp_V$ as asserted. Then (iii) implies that $h\vee h^\sharp_V\in \F(X)$ and $h\vee h^\sharp_V=h$ on $V$, giving that $h\in\efloc(X)$. The \lq in particular' follows, since locally bounded functions can be locally minorized by (negative) equilibrium potentials, see below. \vspace{1mm}
	
	(c) This is a standard semigroup argument plus some approximation argument. It suffices to consider bounded $h$: in fact,
	$$
	\mathscr{E}_\lambda(h,f)=\lim_{n\to\infty}\mathscr{E}_\lambda(h\vee (-n),f)
	$$
	for any $f\in \F_c(X)$ and $h\vee(-n)$ is $\lambda$--defective, provided $h$ is. By rescaling, we can and will assume that $h\ge -1$. Let $f\in \F_c(X)$, $f\ge 0$ and pick $V,W$ open, $V\Subset W\Subset X$ such that $\supp(f)\subset V$. Let $e_{W}$ be the $\lambda$--equilibrium potential of the set $\overline{W}$, see \cite{fukushima}, which is
	$\lambda$--excessive. Therefore,
	$$
	h^\flat:= h\vee (- e_W)
	$$
	is $\lambda$--defective, $h=h^\flat$ on $W$ and $h^\flat\in\F(X)$, because $-h^\flat$ is $\lambda$--excessive and $0\le-h^\flat \le e_W\in\F(X)$. It follows that
	\begin{align*}
		\mathscr{E}_\lambda(h,f)&=\mathscr{E}_\lambda(h^\flat,f)\\
		&=-\lim_{t\searrow 0}\left\langle h^\flat, \frac{1}{t}\left(e^{-\lambda t} P_tf -f\right)\right\rangle\\
		&=-\lim_{t\searrow 0}\left\langle \frac{1}{t}\left(e^{-\lambda t} P_t h^\flat - h^\flat\right), f\right\rangle\\
		&\le 0,
	\end{align*}
completing the proof.
\end{proof}

The following characterization of \emph{$\lambda$-subharmonic functions in the distributional sense} in the Riemannian case served as the blueprint for our new notion of subharmonicity in Definition \ref{defect} below. Recall that a distribution $g$ defined on an open set $U$ of a Riemannian manifold $X$ with Laplace-Beltrami operator $\Delta$ is called \emph{$\lambda$--harmonic}, if $(\Delta  -\lambda) g=0$ on $U$ in the sense of distributions. By local elliptic regularity, this is equivalent to $g$ being a smooth classical solution of the latter equation and thus also to the weak $\lambda$-harmonicity of $g$.

\begin{proposition}\label{aspspo} Assume $X$ is a connected Riemannian $m$-manifold (considered as a smoothing strongly local Dirichlet space in the sense of Example \ref{esxcv}) and let $f\in L^1_\loc(X,\m)$, $\lambda\geq 0$. The following two properties are equivalent: 
\begin{itemize}
\item [(i)] $(\Delta -\lambda)f\geq 0$ in the sense of distributions, that is, $\left\langle  f,   (\Delta-\lambda) \phi \right\rangle \geq 0$ for all $0\leq \phi\in C^\infty_c(X)$,
\item[(ii)] For all $U\Subset X$ there exists a $\lambda$--harmonic function $g\in C_b(U)$ on $U$ such for all $t> 0$ one has $e^{-t\lambda}P^{U}_t (f|_U-g) \geq f|_U-g$.  
\end{itemize}
 \end{proposition}

\begin{proof} Note first that in the above situation, $e^{-t\lambda}P^{U}_t (f|_U-g) \geq f|_U-g$ for all $t> 0$ is equivalent to 
	$$
	t\mapsto e^{-t\lambda}P^{U}_t (f|_U-g)\quad\text{ is increasing in $t> 0$.}
	$$
(i) $\Rightarrow$ (ii): By Theorem D or Theorem 3.2 in \cite{bismar} we can pick open and smooth neighbourhoods $U\Subset U'\Subset  \tilde{U}\Subset X$ and a sequence $(f_k)\in C^{\infty}(\tilde{U})$ such that
\begin{itemize}
	\item $f_k|_{U'}\searrow f|_{U'}$,
	\item $\Delta f_k-\lambda f_k\geq 0$ in $\tilde U$,
		\item $\sup_{U'} f_k|_{U'}\leq \sup_{U'} f|_{U'}+1$ for all $k$.
\end{itemize} 
Note that the first property implies $f_k\to f|_{U'}$ in $L^1(U')$ by dominated convergence and also that $f$ can be chosen upper semicontinuous. By the latter fact we can pick $c\in\IR$ such that $f|_{U'} \leq c/2$ and $f_k|_{U'} \leq c$. Let $\tilde{g}$ be a solution to 
\[
\begin{cases}
(\Delta-\lambda)\tilde{g}=0 &\text{ on }U'\\
\tilde{g}|_{\partial U'}=c
	\end{cases}
\]
and set $h:=f|_{U'}-\tilde{g}$. We have that $h_k:=f_k|_{U'}-\tilde{g}$ is $\lambda$--subharmonic on $U'$ and $h_k\le 0$ on $\partial U'$. Thus $h = \lim_{k\to\infty}h_k\le 0$ on $U'$ by the (weak) maximum principle.
 Moreover, $h_k\to h$ in $L^1(U')$. Let $0\leq \phi\in C^\infty_c(U')$. For all $t>0$ one has
\begin{align*}
(d/dt) \left\langle  e^{-\lambda t}P^{U'}_t h, \phi \right\rangle    &=\left\langle  (\Delta-\lambda) e^{-\lambda t}P^{U'}_t h ,\phi\right\rangle  = \left\langle  e^{-\lambda t}P^{U'}_t h, (\Delta - \lambda) \phi \right\rangle\\
&= \left\langle  h,  e^{-\lambda t}P^{U'}_t\left(\Delta \phi- \lambda \phi\right) \right\rangle = \left\langle  h,  (\Delta - \lambda) e^{-\lambda t}P^{U'}_t\phi\right\rangle.  
\end{align*}
Since $\Delta e^{-\lambda t}P^{U'}_t \phi= e^{-\lambda t}P^{U'}_t \Delta\phi$ is bounded, the latter is
$$
=\lim_k\left\langle  h_k,  \Delta e^{-\lambda t}P^{U'}_t \phi - \lambda e^{-\lambda t}P^{U'}_t\phi\right\rangle.
$$
Next we show
$$
\left\langle  h_k,  \Delta e^{-\lambda t}P^{U'}_t \phi- \lambda e^{-\lambda t}P^{U'}_t\phi \right\rangle\geq 0 \quad\text{for all $k$.}
$$
To this end, we compute
\begin{align*}
	\left\langle  h_k,  \Delta e^{-\lambda t}P^{U'}_t \phi- \lambda e^{-\lambda t}P^{U'}_t\phi \right\rangle&=\int_{\partial U'} h_k \partial_\nu\left( e^{-\lambda t}P^{U'}_t\phi\right)\mathrm{d}\sigma -\int_{U'} (\nabla h_k,\nabla e^{-\lambda t}P^{U'}_t\phi)\mathrm{d}\m\\
	&\phantom{=}-\int_{U'}  \lambda h_k e^{-\lambda t}P^{U'}_t\phi \mathrm{d}\m\\
	& \geq  -\int_{U'} (\nabla h_k,\nabla e^{-\lambda t}P^{U'}_t\phi)\mathrm{d}\m-\int_{U'}  \lambda h_k e^{-\lambda t}P^{U'}_t\phi\mathrm{d}\m\\
	&= \int_{U'} (\Delta  h_k-\lambda h_k)  e^{-\lambda t}P^{U'}_t\phi\mathrm{d}\m-\int_{\partial U'} (e^{-\lambda t}P^{U'}_t\phi) \partial_{\nu} h_k\mathrm{d}\sigma\\
	&\geq 0,
\end{align*}
where we have used $h_k\leq 0$ and $\partial_\nu e^{-\lambda t}P^{U'}_t\phi \leq 0$ (because $e^{-\lambda t}P^{U'}_t\phi \geq 0$) in the first inequality, and $e^{-\lambda t}P^{U'}_t\phi = 0$ on $\partial U'$ together with $(\Delta-\lambda) h_k \geq 0$ and $P^{U'}_t\phi\geq 0$ on $U'$ for the last inequality. We have thus shown $e^{-\lambda t}P^{U'}_t(f|_{U'}-\tilde{g})\geq f|_{U'}-\tilde{g}$ for all $t>0$. Since $P^{U'}_t\to 0$ as $t\to\infty$ strongly in $L^1(U,\m)$ (cf. the proof of Lemma \ref{transfer}(b) below for a detailed argument), we have $f|_{U'}-\tilde{g}\leq 0$. Thus, 
$$
e^{-\lambda t}P^{U}_t(f|_U-\tilde{g}|_U)\geq e^{-\lambda t}P^{U'}_t(f|_{U'}-\tilde{g})\quad\text{on $V$}, 
$$
and the claim follows from setting $g:=\tilde{g}|_U$.\vspace{1mm}

(i) $\Leftarrow$ (ii): Let $U\Subset X$ be arbitrary and pick a $\lambda$--harmonic function $g$ on $U$ as in the assumption. We have $e^{-\lambda t}P^{U}_t(f|_U-g)\to f|_U-g$ in $L^1(U,\m)$, and so $e^{-\lambda t}P^{U}_t(f|_U-g)\to f|_U-g$ as $t\to 0+$ as distributions. Since 
$$
(\Delta-\lambda) e^{-\lambda t}P_t^{U} (f|_U-g)=\partial_t e^{-\lambda t}P^{U}_t(f|_U -g)\geq 0, 
$$
as $e^{-\lambda t}P^{U}_t (f|_U-g)$ is increasing in $t$, we get for all $0\leq\phi\in C^{\infty}_c(U)$,
\begin{align*}	
\int_{U}f(\Delta \phi-\lambda\phi)\mathrm{d}\m&= \int_{U}(f-g)(\Delta \phi-\lambda\phi)\mathrm{d}\m =\lim_{t\to 0+}\left\langle e^{-\lambda t}P^{U}_t(f|_U-g),\Delta \phi-\lambda\phi\right\rangle\\
&= \lim_{t\to 0+}\left\langle (\Delta-\lambda) e^{-\lambda t}P^{U}_t(f|_U-g) ,  \phi\right\rangle\geq 0, 
\end{align*}
and the proof is complete.
\end{proof}

The previous result motivates:

\begin{definition}\label{defect} Let $f\in L^1_\loc(X,\m)$, $\lambda\geq 0$. We say that $f$ is 
\begin{itemize}	
\item \emph{locally $\lambda$--shift defective}, if for all $U\Subset X$ there exists a weakly $\lambda$--harmonic function $g\in  C_b(U)$ such that 
\begin{align}\label{slo} 
e^{-t\lambda}P^{U}_t (f|_U-g) \geq f|_U-g\quad\text{for all $t>0$},
\end{align}
\item \emph{locally $\lambda$--shift excessive}, if $-f$ is locally $\lambda$--shift defective. 
\end{itemize}
\end{definition}

If $\lambda=0$, then we will simply talk about \lq locally shift defective (excessive)'. \vspace{1mm}

\begin{remark} 1. We stress that no a priori regularity on $f$ is needed in the above definition, other than the somewhat minimal $L^1_\loc$-assumption.

2. Note that in the definition of local $\lambda$--shift defectiveness we do not require explicitly that $f|_U-g\leq 0$. The latter will however follow automatically in the locally smoothing case (cf. the transfer principle below). In other words, in the locally smoothing case, $f\in L^1_\loc(X,\m)$ is locally $\lambda$--shift defective, if and only if for all $U\Subset X$ there exists a $\lambda$--harmonic function $g\in  C_b(U)$ such that $f|_U-g$ is $\lambda$--defective in $U$, hence $f|_U-g\leq 0$. This justifies the name \lq locally $\lambda$--shift defective'. \vspace{1mm}
	
3. A function $f\in L^1_\loc(X,\m)$ is obviously locally $\lambda$--shift defective, if and only if $f|_U$ \emph{is locally $\lambda$--shift defective for every $U\Subset X$,} that is, with respect to $(U,\m,\mathscr{E},\mathscr{F}(U))$. Moreover, the class of locally $\lambda$--shift defective functions is stable under taking sums, multiplication by positive constants, and shifting by negative constants. 
\end{remark}

For locally smoothing spaces we get the following result, which will allow us to transfer results for signed functions to functions of varying sign:

\begin{lemma}\label{transfer}\emph{(Transfer principle)} Let $(X,\m,\mathscr{E},\mathscr{F}(X))$ be locally smoothing and  $\lambda\ge 0$.\\
	\emph{(a)} For every $c\in\IR$, $U\Subset X$ there exists $g\in C_b(U)$ weakly $\lambda$-harmonic on $U$, such that $g\geq c$.\\
	\emph{(b)}	If $f\in L^1_\loc(X,\m)$ is locally $\lambda$--shift defective, then with $U$, $g$ as in Definition \ref{defect} one has $f|_U-g\leq 0$.
\end{lemma}

\begin{proof} (a) Pick an open set $U'\subset X$ such that $U\Subset U'$ and $U'\setminus \overline{U}\neq \emptyset$, and such that $(U',\m,\mathscr{E},\mathscr{F}(U'))$ is doubly Feller and irreducible. Pick $0\leq \psi\in C_c(U')\setminus\{0\}$ with $\mathrm{supp}(\psi)\subset U'\setminus \overline{U}$. \\
Assume first $\lambda>0$. Then we can define 
	$$
	0<h:=(H^{U'}+\lambda)^{-1}\psi\in C_0(U')\cap \dom(H^{U'}).
	$$
	Since $h$ is continuous and the closure of $U$ is a compact subset of $U'$, we have $\inf h|_{U} >0$, and so for some $a>0$ we have $c\leq g:= a\cdot h|_U\in C_b(U)\cap \mathscr{F}_\loc(U)$, and given $\phi \in \mathscr{F}_c(U)$ we have 
	\begin{align*}
	\mathscr{E}(g,\phi)+\lambda \left\langle g,\phi\right\rangle = a\mathscr{E}(h,\phi)+\lambda a\left\langle h,\phi\right\rangle =a\langle H^{U'}h,\phi\rangle+a\lambda \left\langle h,\phi\right\rangle=a\left\langle \psi,\phi\right\rangle=0.
	\end{align*}
For  $\lambda =0$, we can simply take a constant function.

(b) For all $t\ge 1$ we have 
$$
e^{-t\lambda}P^{U}_{t-1} P^{U}_1 (f|_U-g) \geq f|_U-g,
$$
where $P^{U}_1 (f|_U-g)\in L^2(U,\m)$. By spectral calculus we have that $P^{U}_{s}$ converges as $s\to\infty$ strongly in $L^2(U,\m)$ to the projection onto $\mathrm{Ker}(H^U)$, the latter space being equal to $\{0\}$, as $H^U$ has a spectral gap. 	
\end{proof}

We continue with:

\begin{lemma}\emph{(Approximation lemma)}
	\label{approx} Assume $(X,\m,\mathscr{E}, \mathscr{F}(X))$ is locally smoothing and $\lambda\geq 0$. Let $f$ be locally $\lambda$--shift defective. \\
	Given $U\Subset X$, pick a weakly $\lambda$--harmonic function $g\in C_b(U)$ on $U$ as in Definition \ref{defect}. Then, for every $k\in\IN$,
	$$
	f^U_k:= e^{-\lambda/k}P^{U}_{1/k}(f|_U-g)+g\in   \F_\loc(U)\cap C_b(U)
	$$
	is locally $\lambda$--shift defective and weakly $\lambda$--subharmonic, and one has $f^U_k\downarrow f|_U$ $\m$-a.e., and so $\pt{f}_U:=\inf_{k\in\IN}f^U_k$ is an upper semicontinuous (u.s.c.) $\m$-representative of $f|_U$. \\
	In particular, the function $f$ has an u.s.c. $\m$-representative.  
\end{lemma}

\begin{proof} Let $V\Subset U$. We have
	\begin{align*}
		&e^{-\lambda t}P^{V}_{t}(f^U_k|_V-g|_V)=e^{-\lambda t}P^{V}_{t}\Big(e^{-\lambda/k}(P^{U}_{1/k}(f|_U-g))|_V \Big)\geq \Big( e^{-\lambda t}P^{U}_{t}(e^{-\lambda/k}P^{U}_{1/k}(f|_U-g))\Big)|_V\\
		&=\Big(e^{-\lambda (t+1/k)}P^{U}_{t+1/k}(f|_U-g)\Big)|_V\geq \Big(e^{-\lambda/k}P^{U}_{1/k}(f|_U-g)\Big)|_V= f^U_k|_V-g|_V,
	\end{align*}
	where have used that $f|_U-g\leq 0$ for the first estimate (part (b) of transfer principle), and that $f$ is locally $\lambda$--shift defective. Thus $f^U_k$ is locally $\lambda$--shift defective. Since $f^U_k-g\le 0$ is continuous and $\lambda$--defective, it is weakly $\lambda$--subharmonic by Theorem \ref{esqq0} (c). This sequence converges pointwise in a monotonically decreasing way on $U$, and the limit function is u.s.c. and coincides $\m$-a.e. with $f|_U$.\vspace{1mm}
	
	The 'in particular' follows with a partition of unity argument. 
\end{proof}

Using the transfer principle, we can establish the following results analogous to those of Theorem \ref{esqq0}, now for functions of varying sign: 

\begin{theorem}\label{esqq} Let $\lambda \geq 0$ and let $(X,\m,\mathscr{E}, \mathscr{F}(X))$ be locally smoothing.\vspace{1mm}

\emph{(a)} If $f\in L^1_\loc(X,\m)$ is locally $\lambda$--shift defective and in $\mathscr{F}_\loc(X)$, then $f$ is weakly $\lambda$--subharmonic. In particular, if $f\in L^\infty_\loc(X,\m)$ is locally $\lambda$--shift defective, then $f$ is weakly $\lambda$--subharmonic \vspace{1mm}
	
\emph{(b)} If $f\in\mathscr{F}_\loc(X)$ is weakly $\lambda$--subharmonic and locally bounded above, then $f$ is locally $\lambda$--shift defective. 
\end{theorem}
\begin{proof} (a) Let $U\Subset X$ be arbitrary and pick $g$ as in Definition \ref{defect}. Then $f|_U-g\in \mathscr{F}_\loc(U)$ is $\lambda$--defective on $U$ by the transfer principle; applying Theorem \ref{esqq0}(c) with $h=f|_U-g$ on $U$ (instead of $X$) shows that this function is weakly $\lambda$--subharmonic on $U$, and so $f|_U=f|_U-g+g$ is weakly $\lambda$--subharmonic on $U$.\vspace{1mm}
	
(b) Let $U\Subset X$ be arbitrary and pick $c\in\IR$ with $f|_U\leq c$ and using the transfer principle pick $g\in C_b(U)$ weakly $\lambda$--harmonic on $U$  with $g\geq c$. Then $f|_U-g\leq 0$ is weakly $\lambda$--subharmonic on $U$, and so applying Theorem \ref{esqq0} (a) to $h=f|_U-g$ on $U$ gives that 
$$
e^{-\lambda t}P^U_t(f|_U-g)\geq f|_U-g,
$$
completing the proof.
\end{proof}

\begin{remark}\label{bemi}The local boundedness from above, required in (b) of the previous Theorem, can be typically deduced from elliptic or parabolic subsolution estimates, which are valid under some mild local assumptions on the space. For the most general result to this effect that we are aware of, we can refer to Theorem 4.7 in \cite{lierl}: under the assumption that the space $(X,\m, \mathscr{E}, \mathscr{F}(X))$ has the property that on every $U\Subset X$ a volume doubling up to a finite scale and a Poincare inequality up to a finite scale are satisfied, a parabolic subsolution estimate is valid, which shows that every weakly $\lambda$-subharmonic function is automatically locally bounded above. Indeed, if $f$ is weakly $\lambda$-subharmonic, then so if $f_+$ by Lemma \ref{sswwq}, and then applying Theorem 4.7 in \cite{lierl} to $(t,x)\mapsto e^{-t\lambda} f_+(x)$, we get that on such a space $f_+$ is automatically locally bounded, and so $f$ is locally bounded above. This applies to all examples of locally smoothing spaces above.
\end{remark}

\section{Kato-Brezis theorem and regularity for $\lambda$-shift defective functions}\label{sec4}

The starting point of this section is the following probabilistic characterization of locally $\lambda$--shift defective functions:

\begin{proposition}\label{prob} Assume that $(X,\m,\mathscr{E}, \mathscr{F}(X))$ is locally smoothing and that $\lambda\geq 0$. A function $f\in L^1_\loc(X,\m)$ is locally $\lambda$--shift defective, if and only if for all $U\Subset X$ there exists a u.s.c $\m$-representative $\pt{f}_U$ of $f|_U$ such that for all $V\Subset U$ with $U\setminus\overline{V}\neq \emptyset$ one has
\begin{align}\label{kaz}
	\mathbb{E}^x_U\left[e^{-\lambda \tau_V}\pt{f}_U(\mathbb{X}_{\tau_V})  \right]\geq \pt{f}_U(x)\mbox{  for q.e. }x\in V.
\end{align}
\end{proposition}

The proof is based on the following technical result by Chen and Kuwae, which holds on any strongly local Dirichlet space and which follows from applying Theorem 2.9 in \cite{kuwae} (see also \cite{chen}) to the Dirichlet form which is shifted by $\lambda$, noting that the result therein holds for nonlocal Dirichlet forms:
	
\begin{lemma}\label{chen} A function $h\in L^\infty_\loc(X)$ is weakly $\lambda$-subharmonic, if and only if for all $V\Subset X$ with $X\setminus\overline{V}\neq \emptyset$ one has
	\begin{align}\label{kua}
		\mathbb{E}^x\left[e^{-\lambda \tau_V}h(\mathbb{X}_{\tau_V})  \right]\geq  h(x)\mbox{  for q.e. }x\in V.
	\end{align}
\end{lemma}

Note that, if $h\in L^\infty_\loc(X)$ is weakly $\lambda$-harmonic, one gets an equality in (\ref{kua}). 

\begin{proof}[Proof of Proposition \ref{prob}] $\Rightarrow$: Let $U\Subset X$. Consider the functions $g$, $f_k^U$, $\pt{f}_U$ from the approximation Lemma \ref{approx}. For every $k\in\IN$ the function $f_k^U-g\le 0$ is continuous and $\lambda$--defective on $U$, in particular, weakly $\lambda$-subharmonic on $U$ by Theorem \ref{esqq0}(c). It follows from using the previous lemma with $h=f_k^U-g$ that
\begin{align*}
	\mathbb{E}^x_U\left[e^{-\lambda \tau_V}(f_k^U-g)(\mathbb{X}_{\tau_V})  \right]\geq (f_k^U-g)(x)\mbox{ for q.e. }x\in V.
\end{align*}
For $k\to\infty$ we arrive at
\begin{align*}
	\mathbb{E}^x_U\left[e^{-\lambda \tau_V}(\pt{f}_U-g)(\mathbb{X}_{\tau_V})  \right]\geq (\pt{f}_U-g)(x)\mbox{ for q.e. }x\in V,
\end{align*}
which using the previous lemma with $h=g$ gives 
\begin{align}\label{ddoo}
	\mathbb{E}^x_U\left[e^{-\lambda \tau_V}g(\mathbb{X}_{\tau_V})  \right]=  g(x)\mbox{  for q.e. }x\in V,
\end{align}
proving \eqref{kaz}.\\
$\Leftarrow$: The function $\pt{f}|_U$ is bounded above for given $U\Subset X$. As above, we can pick a weakly $\lambda$--harmonic function $g\in C_b(U)$ such that $\pt{f}|_U-g\leq 0$ on $U$. Given $V\Subset U$ with $U\setminus\overline{V}\neq \emptyset$, the assumption on $f$ together with (\ref{ddoo}) implies 
$$
\mathbb{E}^x_U\left[e^{-\lambda \tau_V}(\pt{f}_U-g)(\mathbb{X}_{\tau_V})  \right]\geq \pt{f}_U(x)-g(x),\mbox{  for q.e. }x\in V.
$$
For arbitrary $n\in\IN$, this implies
$$
\mathbb{E}^x_U\left[e^{-\lambda \tau_V}\left((\pt{f}_U-g)\vee (-n)\right)(\mathbb{X}_{\tau_V})  \right]\geq (\pt{f}_U(x)-g(x))\vee (-n),\mbox{  for q.e. }x\in V.
$$
and so the locally bounded function $(\pt{f}_U-g)\vee(-n)$ is weakly $\lambda$--subharmonic on $U$ by the previous lemma. From Theorem \ref{esqq0}(a) we get that $(\pt{f}_U-g)\vee(-n)$ is $\lambda$--defective on $U$. Since this is true for arbitrary $n\in\IN$, we get that $\pt{f}_U-g$ is $\lambda$--defective on $U$.
\end{proof}

As noted above, Lemma \ref{chen} is a consequence of Theorem 2.9 in \cite{kuwae}. The latter gives a probabilistic characterization of weak $\lambda$-subharmonicity under the a priori assumption that the function is in $L^\infty_\loc$. This assumption is partially used in \cite{kuwae} to verify necessary measurability proporties and to make the expectation on the LHS of (\ref{kua}) well-defined, noting that the proof goes through an intermediate step, giving a characterization in terms of a submartingale property. \\
Note that we do not make an $ L^\infty_\loc$--assumption in Proposition \ref{prob}, replacing the weak $\lambda$--subharmonicity assumption with a local $\lambda$--shift defectiveness assumption in one direction of the equivalence of Lemma \ref{chen}. The main reason for this is that we get a u.s.c. $\m$-respresentative from our local smoothing assumptions on the space. The main strength of Theorem 2.9 in \cite{kuwae}, however, clearly lies in the very remarkable fact that this result applies to nonlocal forms (with some mild additional regularity assumptions).\vspace{1mm}

As an application of Proposition \ref{prob} we immediately get the following generalization of a consequence of Kato's celebrated inequality, Lemma A from \cite{katoin}, which states that for all $f\in L^1_\loc(\IR^n)$ such that $\Delta f\in    L^1_\loc(\IR^n)$, one has
\begin{equation}\label{kato}
\Delta|f|\ge \mathrm{sign}(f)\Delta f
\end{equation}
in the sense of distributions. This readily gives 
\begin{equation}\label{katoplus}
\Delta f_+\ge 1_{\{ f\ge 0\}}\Delta f.
\end{equation}

This was later sharpened and generalized (see Ancona, \cite{ancona} for a thorough discussion) to the case that $\Delta f$ is a Radon measure in the article \cite{brezis} by Brezis and Ponce. It is maybe superfluous to add that some regularity assumption on $\Delta f$ is needed to ensure that the RHS of \eqref{kato} and \eqref{katoplus} are well-defined distributions. An immediate consequence of the latter estimate is that $\Delta f_+$ is a positive distribution, in other words, $f_+$ is subharmonic in the distributional sense. Remarkably, using the results we have obtained so far, we get a considerably strengthened variant of the statement \emph{$f \in L^1_\loc$ is subharmonic in the distributional sense $\Rightarrow$ $f_+$ is subharmonic in the distributional sense} in our general setting, where neither test functions nor distributions are available. Indeed, the the following result follows immediately from Proposition \ref{prob}:

\begin{theorem}\emph{(Kato-Brezis)}\label{brezis} Assume that $(X,\m,\mathscr{E}, \mathscr{F}(X))$ is locally smoothing and that $\lambda\geq 0$. If $f,g\in L^1_\loc(X,\m)$ are locally $\lambda$--shift defective, then so is $f\vee g$. If $0\leq f \in L^1_\loc(X,\m)$ is locally $\lambda$--shift defective, then $f$ is locally shift defective.
\end{theorem}

The reader may find similar statements in \cite{piggio} for Riemannian manifolds and in \cite{meti} for Carnot groups.\vspace{2mm}

As a first application of the above theorem, we are going to prove:

\begin{theorem}\emph{(Regularity theorem)}\label{main} Let $(X,\m,\mathscr{E}, \mathscr{F}(X))$ be locally smoothing, let $\lambda\geq  0$ and let $ f\in L^1_\loc(X,\m)$ be locally $\lambda$-shift defective. Then for every $c\geq 0$, $q>1$ one has $(f-c)_+^{q/2}\in\F_\loc(X)$, and for $q\geq 2$ this function is weakly subharmonic. 
\end{theorem}

\begin{proof} The function $f-c$ is locally $\lambda$-shift defective and so by Kato-Brezis $(f-c)_+$ is locally $\lambda$-shift defective. Thus we can assume $c=0$ and $f\geq 0$. Then, again by Kato-Brezis, $f$ is locally shift defective. As $f$ has a u.s.c. $\m$-representative, $f$ is locally bounded above (and by $f\geq 0$ then locally bounded) and thus weakly subharmonic by Theorem \ref{esqq} (a). It then follows immediately from the chain rule that for $q\geq 2$ the function $f^{q/2}$ is weakly subharmonic as well (cf. the calculation in the proof of Lemma 2 a) in \cite{sturm1}).\\
	For arbitrary $q>1$ we proceed as follows to show that $f^{q/2}\in\F_\loc(X)$: let $V\Subset X$ be arbitrary and pick $U$ with $V\Subset U\Subset X$. We also pick $0\leq \phi\in C_c(U)\cap \mathscr{F}(U)$ such that $1_V\leq \phi\leq 1_U$, and a sequence $\seqx{f^U}{k}$ as in Lemma \ref{approx} and set $0\leq f|_U\leq f_k:=f^U_k\in \mathscr{F}_\loc(U)\cap C_b(U)$ for notational convenience. Replacing $f_k$ with $f_k+1/k$ if necessary, we can assume $f_k>0$ and still get $f_k\downarrow f|_U$. We have $f_k^{q-1}\in  \mathscr{F}_\loc(U)$ by the chain rule, so
	$$
	0\leq \psi_k:= \phi^2 f_k^{q-1}\in \mathscr{F}(U)\cap C_c(U)
	$$
	is a test function. Since $f_k$ is weakly subharmonic, 
	\begin{align*}
		0\geq \mathscr{E}(f_k,\psi_k)&=\int_U \mathrm{d}\Gamma(f_k,\phi^2f_k^{q-1})\\
		&=\int_U \phi^2\mathrm{d}\Gamma(f_k,f_k^{q-1})+\int_U f_k^{q-1}\mathrm{d}\Gamma(f_k,\phi^2)\\
		&=(q-1)\int_U \phi^2f_k^{q-2}\mathrm{d}\Gamma(f_k)+2\int_U f_k^{q-1}\phi\mathrm{d}\Gamma(f_k,\phi),
	\end{align*}
	where have used the chain rule. Using the Cauchy-Schwarz inequality (\ref{CS}) with
	$$
	v=f_k, \quad w=\phi,\quad f=\phi f_k^{q/2-1},\quad g= f_k^{q/2},
	$$
 on the second term	we get, for $\epsilon>0$:
 \begin{align*}
     2\left|\int_U f_k^{q-1}\phi\mathrm{d}\Gamma(f_k,\phi)\right|&\le 
     2 \left(\int_U \phi^2 f_k^{q-2}\mathrm{d}\Gamma(f_k)\right)^\frac12
     \left(\int_U f_k^{q}\mathrm{d}\Gamma(\phi)\right)^\frac12\\
     &\le \epsilon \int_U \phi^2 f_k^{q-2}\mathrm{d}\Gamma(f_k) +
     \frac{1}{\epsilon}\int_U f_k^{q}\mathrm{d}\Gamma(\phi)
 \end{align*}
 For  $\epsilon>0$
 small enough (depending only on $q$), we deduce that 
	$$
	0\geq (q-1-\epsilon) \int_U \phi^2f_k^{q-2}\mathrm{d}\Gamma(f_k)-\epsilon^{-1} \int_U f_k^{q}\mathrm{d}\Gamma(\phi),
	$$
	Using the chain rule
	$$
	\mathrm{d}\Gamma(f_k^{q/2})=\frac{q^2}{4} f_k^{q-2} \mathrm{d}\Gamma(f_k),
	$$
	we end up with the second estimate in 
	\begin{align*}
		\int_V \mathrm{d}\Gamma(f_k^{q/2})\leq
		\int_U \phi^2\mathrm{d}\Gamma(f_k^{q/2})\leq \frac{q^2}{4\epsilon(q-1-\epsilon)}\int_U f_k^q\mathrm{d}\Gamma(\phi)\leq \frac{q^2}{\epsilon(q-1-\epsilon)}\int_U f_1^q\mathrm{d}\Gamma(\phi),
	\end{align*}
	the latter expression being finite, as $f_1^q\in L^\infty(U,\m)$. Finally, applying Lemma \ref{peter} with $w_k:=f_k^{q/2}|_V$, we can conclude $f^{q/2}|_V\in \mathscr{F}_\loc(V)$, as clearly $w_k\to f^{q/2}|_V$ in $L^2_\loc(V,\m)$ (in fact, in $L^2(V,\m)$).
\end{proof}

\section{Maximum principle}\label{sec5}

The main result of this section is:

\begin{theorem}[Maximum principle]\label{th:maximumprinciple} Let $(X,\m,\mathscr{E}, \mathscr{F}(X))$ be locally smoothing and let $f\in L^{1}_{\loc}(X,\m)$.\\
	a) If $f$ is locally shift defective, then for every $V \Subset X$ there exists a u.s.c $\m$-representative $\pt{f}_V$ on $\overline{V}$ of $f|_V$ such that
	$$
	\sup_{\overline{V}}\pt{f}_V= \sup_{\partial V}\pt{f}_V.
	$$
	b) If $\lambda \geq 0$ and if $f$ is be locally $\lambda$--shift defective, then for every $V\Subset X$ there exists a u.s.c. $\m$-representative $\pt{f}_{+,V}$ on $\overline{V}$ of $f_+|_V$ such that 
	$$
	\sup_{\overline{V}}\pt{f}_{+,V}= \sup_{\partial V}\pt{f}_{+,V}.
	$$
	\end{theorem}

Note that part b) follows immediately from part a), as by the Kato-Brezis theorem, $f_{+}$ is locally shift defective. Note also that it follows immediately from Proposition \ref{prob} that 
$$
\pt{f}_V\leq \sup_{\partial V}\pt{f}_V\quad\text{q.e. in $V$.}
$$
In this sense, the main point of Theorem \ref{th:maximumprinciple} is to show that one can guarantee a pointwise inequality of this form.

\begin{proof}[Proof of Theorem \ref{th:maximumprinciple} a)] From the definition of local shift defectiveness, it is easy to see that it is enough to prove that for every defective function $h$, which is defined in an open neighbourhood $U\Subset X$ of $V$, there exists an $\m$-representative $\pt{h}$ of $h$ such that
\[
\sup_{\partial V}\pt{h}\ge \sup_V\pt{h}.
	\]
As $P^U_{t}h(x)$ is continuous in $x$ and increasing in $t$, we define
$$
\pt{h}:=\lim_{n\to\infty}P^U_{1/n}h.
$$
Assume, by contradiction, that there exists $\bar x \in V$ such that $\pt{h}(\bar x)>\sup_{\partial V} \pt{h}$. We choose $\delta_{1},\delta_{2}\in\IR$ such that 
	\[
	\pt{h}(\bar x)>\delta_1 >\delta_2 > \sup_{\partial V} \pt{h}.
	\]
	Since $\pt{h}$ is u.s.c., we have that $\pt{h}<\delta_2$ in some open neighborhood $N_{\partial V}$ of $\partial V$ in $\bar V$ and, up to choosing a smaller set, we can always assume that 
	\[
	\pt{h}\le\delta_2 \text{ in }\widebar{ N_{\partial V}}.
	\]
	We claim that, by a variation of the proof of Dini's theorem, there exists $t_0>0$ small enough such that,  for all $0<t<t_0$,
	\[
	P_t^{U} h < \delta_1 \text{ in }\widebar{N_{\partial V}}.
	\]
	Indeed, let 
	$$
	E_t:= \{x \in \bar V: P_t^{U} h <\delta_1\}.
	$$ 
	Since $h$ is defective, $E_t\subseteq E_s$ for $t>s$, and by pointwise convergence we get
	$$
	\widebar{N_{\partial V}}\subset \bigcup_{t>0}E_t.
	$$
	Hence $E_{t_0}\supset \widebar{N_{\partial V}}$ for some $t_0$ small enough, and the claimed property follows.\\
	On the other hand
	\[
	P_t^{U} h(\bar x)\ge \pt{h}(\bar x)>\delta_1.
	\]
	Now, $P_t^{U}h$ is defective on $U$ and in $\mathscr{F}(U)\cap C_b(U)$. Thus, Theorem \ref{esqq0} (c) implies that $P_{t}^{U}h$, and hence $u: =  P_t^{U}h - \delta_{1}$, is weakly subharmonic. Since $u \leq 0$ in $\bar V \setminus Y$ where $Y:=\bar V \setminus {N_{\partial V}}$, the weak maximum principle of \cite[Lemma 4]{sturm1} gives that $u\leq 0$ on $V$. In particular,
	\[
	P_{t}^{U} h(\bar x) \leq \delta_{1}.
	\]
	The contradiction completes the proof.
	
\end{proof}

\begin{example}\label{rcdex} Assume that $(X,\m,\mathscr{E}, \mathscr{F}(X))$ stems from an $\mathrm{RCD}^*(0,N)$ space for some $N\geq 1$, see Example \ref{rcd} above. Since $|\nabla \cdot |_*$ is local, one can define $|\nabla f|_*$ also for $f\in \F_\loc(X)$. Let $f$ be weakly $\lambda$-harmonic for some $\lambda\geq 0$, with $\lambda f\in L^q_\loc(X,\m)$ for some $q\geq \max(N,4)$. It then follows from the localized Bochner inequality (cf. \cite{local}) that $|\nabla f|^2_*$ is weakly subharmonic. By Remark \ref{bemi} this function is locally bounded (above), and so Theorem \ref{esqq} (b) shows that this function is locally shift defective. We get that for every $V\Subset X$ there exists a u.s.c. $\m$--representative $\pt{|\nabla f|^2_*}_V$ of $|\nabla f|^2_*|_V$ such that 
$$
\sup_{\overline{V}}\pt{|\nabla f|^2_*}_V = \sup_{\partial V} \pt{|\nabla f|^2_*}_V.
$$
\end{example}

\section{A proof of the BMS-conjecture on locally smoothing spaces}\label{sec6}

We still assume that $(X,\m,\mathscr{E}, \mathscr{F}(X))$ is a strongly local Dirichlet space. Then there exists a possibly degenerate pseudo metric $\dm^\mathscr{E}$ on $X$ (that is, it may happen that $\dm^\mathscr{E}(x,y)=\infty$ for some $x,y\in X$, or that $\dm^\mathscr{E}(x,y)=0$ for some $x,y\in X$ with $x\neq y$), called the \emph{intrinsic metric}, given by 
$$
\dm^\mathscr{E}(x,y):=\sup\{f(x)-f(y):f\in \mathscr{F}_\loc(X)\cap C(X), \mathrm{d}\Gamma(f)\leq \mathrm{d}\m\},
$$
where the condition $\mathrm{d}\Gamma(f)\leq \mathrm{d}\m$ means that $\Gamma(f)$ is absolutely continuous with respect to $\m$, with a density $\leq 1$.  

\begin{definition} We say that $(X,\m,\mathscr{E}, \mathscr{F}(X))$ is
\begin{itemize}
	\item \emph{strictly local}, if $\dm^\mathscr{E}$ induces the original topology on $X$,
	\item \emph{intrinsically complete}, if this space is strictly local and complete. 
\end{itemize}	
\end{definition}
	
If $(X,\m,\mathscr{E}, \mathscr{F}(X))$ is strictly local, then intrinsic completeness is equivalent to properness. This follows from the Hopf-Rinow theorem for locally compact length spaces, in combination with the fact that strictly local Dirichlet spaces are lenght spaces \cite{peter}. 

In spirit of \cite{sturm1}, we say that $(X,\m,\mathscr{E}, \mathscr{F}(X))$ has the \emph{$L^q$-Liouville property}, for given $q\in [1,\infty]$, if every weakly subharmonic function $0\leq f\in L^q(X,\m)$ is constant. One of the main results from Sturm's seminal paper \cite{sturm1} states:

\begin{theorem}\label{leo} If $(X,\m,\mathscr{E}, \mathscr{F}(X))$ is intrinsically complete and irreducible, then it has the $L^q$-Liouville property for every $q\in (1,\infty)$.
\end{theorem}

The following definition is adapted from \cite{batu2}, where the setting of Riemannian manifolds has been considered:

\begin{definition} Given $q\in [1,\infty]$, one says that $(X,\m,\mathscr{E}, \mathscr{F}(X))$ is \emph{$L^q$-positivity preserving}, if for every locally $1$-shift-excessive function $f\in L^q(X,\m)$ one has $f\geq 0$.
\end{definition}

With this definition, we easily get:

\begin{proposition}\label{poss} If $(X,\m,\mathscr{E}, \mathscr{F}(X))$ is locally smoothing with the $L^q$-Liouville property for some $q\in (1,\infty)$, then it is $L^q$-positivity preserving.
\end{proposition}

\begin{proof} If $f\in L^q(X,\m)$ is locally $1$-shift excessive, then $-f$ is locally $1$-shift defective, so $0\leq (-f)_+\in L^q(X,\m)$ is weakly subharmonic by the regularity Theorem \ref{main}; thus this function must be constant by the $L^q$-Liouville property, and thus $=0$.
\end{proof}

Finally, as stated in the introduction, Theorem \ref{leo} in combination with Proposition \ref{poss} immediately leads to a proof of a conjecture by Braverman, Milatovic and Shubin extended to our general setting (see Conjecture P in Appendix B, p. 679, and also the discussion in \cite{batu3}):

\begin{theorem}\emph{(BMS-conjecture)} If $(X,\m,\mathscr{E}, \mathscr{F}(X))$ is locally smoothing, intrinsically complete and irreducible, then it is $L^q$-positivity preserving for all $q\in (1,\infty)$.
\end{theorem}

Note that every connected and complete Riemannian manifold, every RCD$^*$--space, as well as every Carnot group satisfies the assumptions of this theorem. Indeed, in the Riemannian case the intrinsic completeness is well-known, in the RCD$^*$--case it is shown in \cite{ags}, and for the Carnot group in \cite{dragoni}. In all cases, the irreducibility follows from the strict positivity of the heat kernel.

%
\setcounter{section}{0}
\renewcommand*{\thesection}{\Alph{section}}

\section{Appendix: A class of locally smoothing spaces}

We assume in this section that $(X,\m,\mathscr{E}, \mathscr{F}(X))$ is a strongly local Dirichlet space, with $(X,\dm)$ a metric space such that for all open balls $\B(x,r)$ one has
$$
\m(x,r):=\m(\B(x,r))<\infty.
$$

We assume furthermore that $t\mapsto P_t$ has a strictly positive Markovian heat kernel, by which we mean that there exists a jointly continuous function
$$
p:(0,\infty)\times X\times X\longrightarrow (0,\infty)
$$
with 
$$
\int_Xp(t,x,y)\Id\m (y)= 1,\quad p(t+s,x,y)=\int_X p(t,x,z)p(s,z,y) \mathrm{d}\m(z)\quad\text{for all $t,s>0$, $x,y\in X$,}
$$
such that for all $t>0$, $0\leq f\in L^2(X,\m)$ one has
$$
P_tf(x)=\int_X p(t,x,y) f(y)\Id\m(y)\quad\text{for $\m$-a.e. $x\in X$.}
$$
The same formula then gives a pointwise well-defined representative of $P_tf$ for all $f\in L^q(X,\m)$, $q\in [1,\infty]$, or for nonnegative or nonpositive $f\in \mathscr{B}(X,m)$.\\
The main result of this section is:

\begin{proposition}\label{app} Assume that in the above situation the following two assumptions are satisfied: 
\begin{itemize}
\item \emph{(Local volume doubling)} There exists constants $C>0$, $\alpha\geq 1$, such that for all $0<r<R$, $x\in X$ one has 
$$
\frac{\m(x,R)}{\m(x,r)}\leq C e^{CR} (R/r)^\alpha.
$$
\item \emph{(Local upper Li-Yau heat kernel bound)} There exist constants $c_{1},c_2,c_{3},>0$, such that for all $t>0$, $x,y\in X$ one has 
$$
p(t,x,y)\leq c_{1} \m(x,\sqrt{t})^{-1}e^{-c_2\frac{\dm(x,y)^2}{t}}e^{c_3t}.
$$
\end{itemize}
Then for all $U\Subset X$ the space $(U,\m,\En,\F(U))$ is ultracontractive, and $(X,\m,\mathscr{E}, \mathscr{F}(X))$ is doubly Feller and irreducible. In particular, $(X,\m,\mathscr{E}, \mathscr{F}(X))$ is locally smoothing.
\end{proposition}	

\begin{proof} Let us first record the following consequences of local volume doubling: firstly, there exists $D>0$ such that for all $x,y\in X$, and all $t>0$ one has
	\begin{align}\label{doub}
		\frac{\m(x,\sqrt{t})}{\m(y,\sqrt{t})}\leq D e^{Dt}e^{\frac{\dm(x,y)^2}{t}}.
	\end{align}
This fact is wellknown (cf. Section 2 of \cite{baum} for the simple argument). \\
Secondly, there exist $A,B>0$ such that for all $\delta,t\in [0,1)$, $x\in X$, one has 
	\begin{align}\label{sa}
	\int_{\{\dm(x,y)\geq \delta\}}e^{-\frac{-\dm(x,y)^2}{ct}}\Id\m(y)\leq \m(x,\sqrt{t})A e^{-\frac{\delta^2}{Bt}}.
	\end{align}
This is consequence of Lemma 5.2.13 in \cite{saloff-buch} and local volume doubling. \vspace{1mm}

1. Let $U\Subset X$. The space $(U,\m,\En,\F(U))$ is ultracontractive: indeed, given $0\leq f\in L^1(U,\m)$ one has
$$
P^U_tf(x)\leq P_tf(x)= \int_U p(t,x,y) f(y) \mathrm{d}\m(y)\leq \Big(\sup_{x',y'\in U}p_t(x',y')\Big) \int_Uf(y)\Id\m(y),
$$
as $(x,y)\mapsto p(t,x,y)$ is continuous.\vspace{1mm}

2. The space $(X,\m,\mathscr{E}, \mathscr{F}(X))$ is doubly Feller: we first show that $P_t:C_0(X)\to C_0(X)$. Clearly, it suffices to show that $P_tf\in C_0(X)$ for all $f\in C_c(X)$. With $K$ the support of $f$, the local Li-Yau bound implies
$$
|P_tf(x)|\leq c_1 e^{c_3t} \sup_{y\in K }\m(\B(y,\sqrt{t}))^{-1}\left\|f\right\|_{\infty}\int_{K}  e^{-c_2\frac{\dm(x,y)^2}{t}} \Id\m(y),
$$
which tends to $0$ as $x\to\infty$ and we can see the continuity of $x\mapsto P_tf(x)$ by dominated convergence. Note that
$$
\sup_{y\in K }\m(\B(y,\sqrt{t}))^{-1}<\infty
$$
is a consequence of (\ref{doub}) and that $\m$ has a full support.\\
Next, we show that $P_t:L^\infty(X,\m)\to C(X)$: since $P_t$ has a jointly continuous integral kernel, Corollary 2.2 in \cite{schilling} (which is formulated for $\IR^n$, but the proofs apply in our situation), show that it suffices to prove that $P_t:C_b(X)\to C(X)$. This again follows (cf. the remark prior to Corollary 2.2 in \cite{schilling}), from the fact that $P_t1$ is continuous, as 
$$
P_t1(x)=\int_Xp(t,x,y)\mathrm{d}\m(y)=1,
$$
and the fact that $P_tf\in C(X)$ for all $f\in C_c(X)$ (which we have already shown above). \\
Next we show that
\begin{align}\label{conn}
	\lim_{t\to 0+}P_tf=f\quad \text{uniformly in every compact $K\subset X$, for all $f\in C_0(X)$.}
\end{align}
To this end, given $\epsilon >0$ we can pick $\delta>0$ such that for all $x\in K$, all $y\in B(x,\delta)$ one has $|f(x)-f(y)|\leq \epsilon$. For all $t>0$ we have
$$
|P_tf(x)-f(x)|\leq \int_{\{\dm(x,y)<\delta\}}p(t,x,y)|f(x)-f(y)|\Id\m(y)+\int_{\{\dm(x,y)\geq \delta\}}p(t,x,y)|f(x)-f(y)|\Id\m(y).
$$
The first summand is
$$
\leq \epsilon\int_{X}p(t,x,y)\Id\m(y)= \epsilon.
$$
The second is estimated using the local Li-Yau bound according to
$$
\leq 2c_{1}e^{c_3t}\left\|f\right\|_\infty \m(B(x,\sqrt{t}))^{-1}\int_{\{\dm(x,y)\geq \delta\}}e^{-c_2\frac{\dm(x,y)^2}{t}}\Id\m(y),
$$
which by (\ref{sa}) for all $0<t<1$ is 
$$
\leq 2c_{1}e^{c_3t}\left\|f\right\|_\infty A e^{-\frac{\delta^2}{Bt}}.
$$
This completes the proof of the doubly Feller property.\vspace{1mm}

3. The space $(X,\m,\mathscr{E}, \mathscr{F}(X))$ is irreducible: indeed, given $0\leq f\in L^2(X)\setminus\{0\}$, this follows from 
$$
P_tf(x)=\int_{X} p(t,x,y)f(y)\mathrm{d}\m(y)>0,
$$
as $(x,y)\mapsto p(t,x,y)>0$ is continuous and $\m$ has a full support.

\end{proof}

We record the following corollary to the proof of the above result:

\begin{proposition}\label{app2} Assume that in the above situation one has 
	\begin{itemize}
		\item $P_t:C_b(X)\to C(X)$ for all $t>0$,
		\item $P_t:C_c(X)\to  C_0(X)$ for all $t>0$,
		\item $\lim_{t\to 0+}P_tf(x)=f(x)$ for all $f\in C_0(X)$, $x\in X$.
		\ 
	\end{itemize}
	Then for all $U\Subset X$ the space $(U,\m,\En,\F(U))$ is ultracontractive, and $(X,\m,\mathscr{E}, \mathscr{F}(X))$ is doubly Feller and irreducible. In particular, $(X,\m,\mathscr{E}, \mathscr{F}(X))$ is locally smoothing.
\end{proposition}

\end{document}